\documentclass[11pt,twoside%,draft
]{article}
\usepackage{amsfonts,amssymb}
\RequirePackage[OT1]{fontenc}
\RequirePackage[aap,amsthm,amsmath,natbib,noinfoline]{imsart}
\RequirePackage[dvips%,colorlinks
]{hyperref}

\RequirePackage{hypernat}
\startlocaldefs
\def\@journal{\mbox{}}

\newcommand{\Prob}{\mathop{{}\mathbb P}{}}
\newcommand{\EE}{\mathop{{}\mathbb E}{}}
\newcommand{\Var}{\mathop{\rm Var}{}}
\newcommand{\NN}{{\mathbb N}}
\newcommand{\ZZ}{{\mathbb Z}}
\newcommand{\RR}{{\mathbb R}}
\newcommand{\RT}{\textsf{RT}}

\newcommand{\diff}{{\rm d}}
\newcommand{\e}{{\rm e}}

\newcommand{\myp}{\mbox{$\:\!$}}
\newcommand{\mypp}{\mbox{$\;\!$}}

\theoremstyle{plain}

\newtheorem{theorem}{Theorem}[section]

\newtheorem{corollary}[theorem]{Corollary}
\newtheorem{lemma}[theorem]{Lemma}

{\end{trivlist}}

\theoremstyle{remark}
\newtheorem*{remark}{Remark}
\newtheorem{example}[theorem]{Example}

\numberwithin{equation}{section}

\def\thanksmark@fmt#1{\hbox{\footnotesize $^{#1}$}}
\endlocaldefs

\begin{document}

\begin{frontmatter}

\title{On the variance of the number of occupied boxes}
%\thankstext{*}{If you want to make thanks or other footnotes}
\runtitle{Variance of the number of occupied boxes}

%\begin{aug}
\author{\fnms{Leonid V.} \snm{Bogachev},\thanksref{t1}
\ead[label=e1]{bogachev@maths.leeds.ac.uk}%
}%
\address{L.~V.~Bogachev\\
Department of Statistics\\
University of Leeds\\
Woodhouse Lane\\
Leeds LS2\;9JT\\
United Kingdom\\
\printead{e1}}%
\author{\fnms{Alexander V.} \snm{Gnedin}%\corref{}\thanksref{t2}%
\ead[label=e2]{gnedin@math.uu.nl, \ yakubovich@math.uu.nl}}
\address{A.~V.~Gnedin and Yu.~V.~Yakubovich\\
Department of Mathematics\\
Utrecht University\\
Budapestlaan 6, De Uithof\\
PO Box 80.010, 3508 TA Utrecht\\
The Netherlands\\ \printead{e2}}\\
\and
\author{\fnms{Yuri V.} \snm{Yakubovich}\corref{}\thanksref{t3}
%\ead[label=e3]{yakubovich@math.uu.nl}
}
%\address{Yu.~V.~Yakubovich\\
%Department of Mathematics\\
%Utrecht University\\
%Budapestlaan 6, De Uithof\\
%O Box 80.010\\
%3508 TA Utrecht\\
%he Netherlands\\ \printead{e3}}

\thankstext{t1}{\myp{}Supported in part by DFG Grant 436 RUS
113/722 and a WUN GEP visiting grant.}%
%\thankstext{t2}{\mypp{}Corresponding author.}
\thankstext{t3}{\myp{}Supported by NWO Open
Competition Grant 613.000.304.}

\affiliation{University of Leeds, University of Utrecht and
University of Utrecht}

\runauthor{L.~V.\ Bogachev, A.~V.\ Gnedin and Yu.~V.\ Yakubovich}

%\end{aug}
%\maketitle

\begin{abstract}
We consider the occupancy problem where balls are thrown
independently at infinitely many boxes with fixed positive
frequencies. It is well known that the random number of boxes
occupied by the first $n$ balls is asymptotically normal if its
variance $V_n$ tends to infinity. In this work, we mainly focus on
the opposite case where $V_n$ is bounded, and derive a simple
necessary and sufficient condition for convergence of $V_n$ to a
finite limit, thus settling a long-standing question raised by
Karlin in the seminal paper of 1967. One striking consequence of our
result is that the possible limit may only be a positive integer
number. Some new conditions for other types of behavior of the
variance, like boundedness or convergence to infinity, are also
obtained. The proofs are based on the poissonization techniques.
\end{abstract}



\end{frontmatter}

\section{Introduction}\label{sec1}
The classical occupancy problem is one of the cornerstones of
discrete probability, dating back to its early ages (and hence
encountered over and over again by the generations of students
studying elementary probability through the evergreen hits like the
birthday problem, the coupon collector's problem, etc.\
\citep{Feller1, Aldous}).
%\citep[Ch.~II]{Feller1}, \citep[Ch.~E]{A}).
It still attracts lots of research interest, especially in recent
years, mainly due to its numerous applications spreading across the
board, from sampling statistics and quality control to quantum
physics, bioinformatics and computer science. For an introduction to
the field and a survey of the many models and results, see
\citep{Johnson-Kotz, %Kolchin-Chistyakov,
Kolchin-Sevastyanov-Chistyakov, Kotz-Balakrishnan,
Ivanov-Ivchenko-Medvedev, 
Charalambides} and further references to
original work therein.

In this paper, we are concerned with a version of the occupancy
problem in an \emph{infinite urn scheme} (first considered by
Bahadur~\citep{Bahadur} and later on studied by Darling
\citep{Darling} and most systematically by Karlin {\citep{Karlin}),
in which the balls labeled $1,2,\dots$ are thrown independently at
an infinite array of boxes (\emph{urns}) $j=1,2,\dots$, with fixed
probability (\emph{frequency}) $p_j$ of hitting box $j$. The
frequencies $p_j$ are assumed to be strictly positive and satisfying
\begin{equation}\label{add1}
\|p\|:=\sum_{j=1}^\infty p_j=1\myp.
\end{equation}
Without loss of generality, we further assume that the sequence
$(p_j)$ is non-increasing, $p_1\ge p_2\ge\cdots$.

Let $K_n$ be the number of boxes discovered by the first $n$ balls
(i.e., occupied by at least one of the first~$n$ balls). Many other
interpretations of this functional appear in the literature: for
instance, when $(p_j)$ is considered as a probability distribution
on positive integers, $K_n$ is the number of distinct values
occurring among $n$ random values sampled independently from
$(p_j)$. Since there are infinitely many boxes, $K_n$ increases
unboundedly (with probability one) as more balls are thrown, which
also implies (e.g., by Fatou's lemma) that the same is true for the
expected number of occupied boxes, $\EE(K_n)$.
%\uparrow \infty$ as $n\to\infty$.
Moreover, as shown by Karlin \citep[Theorem~8]{Karlin},
$\lim_{n\to\infty}K_n/\EE(K_n)= 1$ with probability one (an earlier
result about convergence in probability was obtained by
Bahadur~\cite{Bahadur}).

The more delicate asymptotic properties of the random variable $K_n$
are largely determined by its variance $V_n:=\Var(K_n)$. It is known
\cite{Karlin,Dutko,Hwang-Janson} that the distribution of $K_n$
converges to a normal distribution provided that $V_n\to\infty$ as
$n\to\infty$. The latter occurs, for instance, when the frequencies
have a power-like decay, $p_j\sim cj^{-\alpha}$ ($j\to \infty$) with
$\alpha>1$ or, more generally, satisfy a condition of regular
variation \citep{Karlin}. (Here and throughout, $c$\/ stands for a
generic positive constant, specific value of which is not
important.)

\subsection{Main result\/\textup{:}
the case of converging variance}\label{sec1.1}
In this paper, we essentially focus on the opposite situation, that
is, when $V_n$ is uniformly bounded (and hence the distribution of
$K_n$ does not converge to normal). In particular, we prove the
following surprising characterization of frequencies $(p_j)$ for
which the variance $V_n$ tends to a finite limit as $n\to\infty$.
\begin{theorem}\label{thm:limit}
A finite limit\/ $v:=\lim_{n\to\infty}V_n$ exists if and only if for
some integer\/ $k\ge1$ the frequencies satisfy the ``lagged ratio''
condition
\begin{equation}\label{pv}
\lim_{j\to\infty}\frac{p_{j+k}}{p_j}=\frac{1}{2}\mypp,
\end{equation}
and in this case the limiting value $v$ coincides with the lag~$k$.
\end{theorem}
The striking consequence of this result is that whenever the finite
limit of the sequence $(V_n)$ exists, it must be a positive integer
number, $v\in\NN$\myp.

The issue of converging variance was first queried in the seminal
paper by Karlin \citep{Karlin}, where in particular he appreciated
as ``formidable if not impossible'' the task to determine the
behavior of the variance $V_n$ without some regularity assumptions.
In particular, adopting the condition of regular variation of the
frequency tail, he came up with a sufficient condition for the
existence of a finite limit of $V_n$ \citep[Theorem~2]{Karlin}. In
fact, as we shall see below (in Section~\ref{sec5}), convergence to
a finite limit, combined with the special dyadic structure of the
counting measure controlling the frequency input, is a regularity
condition in itself, being strong enough to ensure the result of
Theorem~\ref{thm:limit}. (To be more precise, the ``dyadic'' feature
mentioned above, pertains primarily to the \emph{poissonized}
version of the problem, i.e., with randomized number of balls, see
Section \ref{sec2} below).

The prototypical (apparently folklore) instance of frequencies
$(p_j)$ with converging variance $V_n$ is the geometric sequence of
ratio $1/2$ (i.e., $p_j=2^{-j}$), where one can show with some
effort that $V_n\to 1$ as $n\to\infty$ (see \citep{Karlin, Dutko,
Hwang-Janson}). Note that our condition (\ref{pv}) is obviously
satisfied here with $k=1$, hence the result. The mechanism leading
to such a simple answer is due to a resonance of the ratio $q=1/2$
with the intrinsic dyadic structure of the variance, resulting in
massive cancelation of oscillating terms (again, in the poissonized
version, see Example \ref{ex:geometric1/2} below). Recently, such
cancelations have been explained directly for the original model
(i.e., for $V_n$) using sophisticated analytic methods
\citep{Prodinger, Archibald-Knopfmacher-Prodinger}.

It seems to be less well known that for generic geometric
frequencies $p_j=c\mypp q^{-j}$, the (finite) limit of $V_n$ exists
if $q=2^{-1/k}$ ($k\in\NN$), with the limiting value $v=k$ (see
\citep[\S\mypp4, page~15]{Janson}). Again, using Theorem
\ref{thm:limit} one gets this answer immediately, together with the
``only if'' statement; moreover, the same conclusion can be readily
extended to sequences $(p_j)$ from the parametric class $\RT_q$ (see
\citep{Burris, Bell-Burris, Granovsky}), defined by the property
\begin{equation}\label{eq:RT}
\lim_{j\to\infty} \frac{p_{j+1}}{p_j}=q\myp,
\end{equation}
thus asymptotically mimicking the geometric decay. (Some concrete
examples of distributions in the $\RT_q$ class, complementing the
geometric instance, will be given below in Section~\ref{sec1.3}.)
Indeed, in the $\RT_q$ case equation (\ref{pv}) amounts to
$q^k=1/2$, whence $q=2^{-1/k}$. Of course, condition (\ref{eq:RT})
is too restrictive for the criterion (\ref{pv}), as can be seen for
instance by merging $k$ geometric sequences of the same ratio
$q=1/2$ (and normalizing the resulting sequence so as to
satisfy~(\ref{add1})).

The following ``decomposition'' interpretation of Theorem
\ref{thm:limit} clarifies the compound structure of frequency
sequences $(p_j)$ that exhibit convergence of the variance. Observe
that by condition (\ref{pv}), the sequence $(p_j)$ splits in a
disjoint fashion into $k$ non-increasing subsequences
$p_j^{(i)}\!:=p_{\myp i+k(j-1)}$ \,($i=1,\dots,k$), each belonging
to the $\RT_{1/2}$ class:
\begin{equation}\label{eq:decomposition}
(p_j)=\bigsqcup_{\,i=1}^{\,k} \bigl(p_j^{(i)}\bigr)\,:\qquad
\lim_{j\to\infty} \frac{p_{j+1}^{(i)}}{p_j^{(i)}}=\frac{1}{2}\qquad
(i=1,\dots,k).
\end{equation}
Moreover, by the ``if'' part of Theorem \ref{thm:limit}, each of the
$k$ constituent subsequences %$\bigl(p_j^{(i)}\bigr)$
brings a unit contribution to the overall limiting variance $v=k$.

Such a decomposition may be interpreted as splitting the initial
array of boxes $1,2,\dots$ into $k$ infinite sub-arrays
$\{i+k(j-1),\ j=1,2,\dots\}$ ($i=1,\dots,k$), and allocating the
balls to boxes in a two-stage procedure as follows: for each ball, a
destination array is chosen independently with probabilities
$\|p^{(i)}\|$, and the ball is then thrown with the corresponding
(re-scaled) frequencies $p_j^{(i)}/\|p^{(i)}\|$ \,($j=1,2,\dots$).
The additivity of the variance in this procedure, as predicted by
Theorem \ref{thm:limit}, may be somewhat surprising, given the
apparent dependence of the partial occupancy numbers $K_n^{(i)}$
($i=1,\dots,k$). However, additivity becomes quite transparent in
the poissonized setting, where the dependence between boxes is
removed (see a remark in Section~\ref{sec2.2}).

\subsection{Geometric frequencies}\label{sec1.2}
Historically, there has been some confusion about the converging
variance in the geometric model. Controversy started in
\,\citep[Example~6]{Karlin},
%page~385
where Karlin asserted that
his sufficient condition for convergence \citep[Theorem~2]{Karlin}
was satisfied for \emph{every} geometric sequence $p_j=c\mypp q^{j}$
($0<q<1$), with the limiting value given by $v=\log_{1/q}2$. As we
have seen, this is false unless $q$ belongs to the countable set
$\{2^{-1/k},\ k\in\NN\}$. A more careful inspection reveals that
Karlin's condition, if applied accurately, does yield the correct
answer in the geometric case, properly discriminating between
convergence vs.\ divergence! Moreover, we have found out, quite
unexpectedly, that Karlin's condition (decorated in \citep{Karlin}
with some superfluous assumptions and originally conceived as just a
sufficient condition) proves to be \emph{necessary and sufficient},
being equivalent to our own criterion proved in
Lemma~\ref{lm:Cesaro}. We will discuss this link below, in Section
\ref{sec5Karlin}.

That there was something wrong with Example~6 in \citep{Karlin} was
subsequently pointed out by Dutko \citep[page~1258]{Dutko}, who
noticed that $V_n$ is bounded below by a positive constant,
uniformly in $n$ and $q$, hence the limit $v=\log_{1/q}2$ cannot be
valid at least for small values of $q$ (when $\log_{1/q} 2$ gets
arbitrarily close to zero). However, Dutko \citep[page~1258]{Dutko}
apparently claimed that the limit of the variance fails to exist for
\emph{each} $q\ne1/2$, thus missing the other values, $q=2^{-1/k}$,
\,$k>1$. Unfortunately, he gave no details to support such a
conclusion, referring to his unpublished thesis \citep{Dutko-PhD},
which is not easily available.

More recent studies \citep{Archibald-Knopfmacher-Prodinger,
Hitczenko-Louchard, Louchard-Prodinger-Ward, Prodinger} have shed
much light on the geometric model. Hitczenko and Louchard
\citep{Hitczenko-Louchard} (motivated by random compositions of
natural numbers) were apparently first to prove analytically that
$V_n=1+o(1)$ in the geometric case with $q=1/2$, contrary to
``popular belief'' \citep{Prodinger} that persistent oscillations
are ubiquitous in discrete random structures involving geometric
distribution (see, e.g.,
\citep{Sedgewick-Flajolet,
Szpankowski,
Hwang-Janson}). Prodinger
\citep{Prodinger} gave an alternative proof of this asymptotics
(along with a similar result for a particular model of data search
trees called PATRICIA tries), proceeding from the general
``oscillatory'' framework. Recently, Archibald \emph{et al.}
\cite[Theorem~2]{Archibald-Knopfmacher-Prodinger} derived a very
precise asymptotic expansion
\begin{equation}\label{Arch}
V_n=\log_{1/q} 2 + \delta_V(\log_{1/q}n)+o(1)\qquad(n\to\infty)\myp,
\end{equation}
where $\delta_V(x):=\delta_E(x+\log_{1/q} 2)-\delta_E(x)$ with
$\delta_E(\cdot)$ periodic of period $1$ and zero mean (the latter
function emerges in a similar expansion for $\Phi_n$, the expected
value of $K_n$). If $q=1/2$ then $\log_{1/q} 2=1$, and from the
expansion (\ref{Arch}) it is seen that the oscillating term vanishes
due to $1$-periodicity of $\delta_E(\cdot)$, since
$\delta_V(x)=\delta_E(x+1)-\delta_E(x)=0$ (see \cite[Appendix~A,
page~1079]{Archibald-Knopfmacher-Prodinger}. In fact, the same
argument is true for any $q=2^{-1/k}$ ($k\in\NN$), when
$\log_{1/q}2=k$ and hence $\delta_V(x)=\delta_E(x+k)-\delta_E(x)=0$
(see~\citep[\S\mypp4, page~15]{Janson}).

\subsection{Bounded variance and convergence to infinity}\label{sec1.3}
One can also wonder about conditions for other possible types of
behavior of the variance $V_n$. We shall prove the following
criterion of uniform boundedness, again set in terms of the lagged
ratio $p_{j+k}/p_j$ compared to the \emph{upper} threshold $1/2$
[cf.\ (\ref{pv})].

\begin{theorem}\label{thm:bound}
The sequence\/ $(V_n)$ is bounded if and only if there exists a
positive integer $k$ such that the frequencies\/ $(p_j)$ satisfy the
condition
\begin{equation}\label{eq:boundsup}
\limsup_{j\to\infty} \frac{p_{j+k}}{p_j}\le \frac{1}{2}\mypp.
\end{equation}
Moreover\myp\textup{,} if\/ $k$ is the least integer with the
property\/ \textup{(\ref{eq:boundsup})}\textup{,} then $(V_n)$
satisfies a sharp asymptotic bound\/ $\limsup_{n\to\infty} V_n\le
k$\myp.
\end{theorem}

This situation is exemplified by the generic geometric frequencies,
with arbitrary ratio $0<q<1$. Another example is the Poisson
frequencies $p_j=c\mypp \lambda^j/j\myp!$ \,($\lambda>0$), where the
variance $V_n$ is bounded but does not converge: indeed, here
$p_{j+k}/p_j\sim (\lambda/j)^k\to 0$ as $j\to\infty$, hence
(\ref{eq:boundsup}) is fulfilled whereas (\ref{pv}) fails. A larger
class is that of quasi-binomial distributions \cite{Kerov}, given by
$p_j=(c/j\myp!)\prod_{i=0}^{j-1} (\lambda+iq)$ with parameters
$\lambda>0$, $0\le q<1$. (To explain the name, note that
$c^{-1}=(1-q)^{-\lambda/q}-1$ for $q>0$, while for $q=0$ one has, in
a continuous fashion, $c^{-1}=\e^\lambda-1$, thus recovering the
Poisson normalization constant.) Somewhat similar but different
parametric family is given by the negative binomial distribution
$p_j=(c\mypp q^j/j\myp!)\prod_{i=0}^{j-1}
(\lambda+i)=c\mypp\binom{\lambda+j-1}{j}\mypp q^j$, with
$\lambda>0$, $0<q<1$ [here $c^{-1}=(1-q)^{-\lambda}-1$\myp].

Note that all these examples belong to classes $\RT_q$ with $0\le
q<1$. It is possible to construct more general examples using the
``decomposition'' reformulation of Theorem \ref{thm:bound} in the
spirit of (\ref{eq:decomposition}), in that the variance $V_n$ is
uniformly bounded if and only if the sequence $(p_j)$ may be split
in a disjoint fashion into a finite number of subsequences, each of
which satisfies condition (\ref{eq:boundsup}) with $k=1$ (e.g., each
from $\RT_{q_i}$ with $0\le q_i\le 1/2$, \,$i=1,\dots,k$).

We shall also address the classical question of convergence to
infinity and produce new conditions ensuring that $V_n\to\infty$.
Note, however, that in contrast to the convergent or bounded cases,
no necessary and sufficient criteria are available without extra
regularity assumptions.  To illustrate our results in this
direction, let us formulate here two sufficient conditions, the
first of which is set in terms of the lagged ratios $p_{j+k}/p_j$
against the \emph{lower} threshold $1/2$ [cf.\ (\ref{eq:boundsup})],
while the second one is based on the ``tail ratio''
\begin{equation}\label{eq:rho}
\rho_j:=\frac{1}{p_j}\sum_{i>j} p_{i}\myp.
\end{equation}

\begin{theorem}\label{th:->infty}
Suppose that for each integer\/ $k\ge1$\myp\textup{,}
\begin{equation}\label{eq:>=}
\liminf_{j\to\infty}\frac{p_{j+k}}{p_j}\ge \frac{1}{2}\mypp.
\end{equation}
Then it follows that
\begin{equation}\label{eq:pjovertail0}
\lim_{j\to\infty} \rho_j=\infty\myp,
\end{equation}
which in turn implies that\/ $V_n\to\infty$ as\/ $n\to\infty$\myp.
\end{theorem}

Examples to Theorem \ref{th:->infty} are immediately supplied by the
class $\RT_1$, where condition (\ref{eq:>=}) is obviously satisfied
for any $k\ge1$. More complex examples (not in $\RT_1$) will be
constructed in Sections \ref{sec4.1} and~\ref{sec4.3}.

\begin{remark}
The tail ratio (\ref{eq:rho}) can be expressed as
$\rho_j=(1-h_j)/h_j$\myp, where $h_j=p_j\big/\sum_{i\ge j}^\infty
p_{i}$ is the discrete-time hazard rate, a key characteristic in
reliability theory and survival analysis (see, e.g., 
\cite{Barlow-Proschan%,Cox-Oakes
}). The latter quantity also appears
in the extreme value theory in connection with records from discrete
distributions, where it is interpreted as the probability that $j$
is a record value (see, e.g., \citep{Vervaat,Nevzorov}).
%\citep{Shorrock,Vervaat,Nevzorov,Stepanov}).
In the occupancy context, condition (\ref{eq:pjovertail0}) is
related to the ``probability of a tie for first place''
$\Prob\{X_{n,M_n}=1\}$, where $M_n:=\max{}\{j: X_{n,j}\neq 0\}$ is
the largest index among the occupied boxes after $n$ throws. More
specifically, it has been proved \citep{Eisenberg-Stengle-Strang,
Baryshnikov-Eisenberg-Stengle} that condition (\ref{eq:pjovertail0})
is satisfied if and only if
\begin{equation}\label{ptfp}
\Prob\{X_{n,M_n}=1\}\to 1\qquad(n\to\infty)\myp,
\end{equation}
and moreover, if (\ref{eq:pjovertail0}) fails then
$\Prob\{X_{n,M_n}=1\}$ does not converge at all.
%It is also known \citep{Qi} that (\ref{eq:pjovertail0}) is necessary
%and sufficient in order that $\EE(X_{n,M_n})\to1$.
This, combined with Theorem~\ref{sctail}, shows that (\ref{ptfp})
implies both $V_n\to\infty$ and $\Phi_{n,1}\to\infty$, which is a
surprising connection between the behavior in the extreme-value
range and the global characteristics of the sample. These facts
equally apply to the poissonized model.
\end{remark}

\subsection{Outline} The rest of
the paper is organized as follows. Section \ref{sec2} contains
general formulas and introduces the poissonization technique. In
Section~\ref{sec3}, we connect the variance $V_n$ with the mean
number of singletons (i.e., the boxes occupied by exactly one of the
first $n$ balls) and derive useful upper bounds. We also obtain here
a basic integral representation of the poissonized variance $V(t)$
via the Laplace transform of the function $\Delta\nu(x)$, counting
the frequencies $p_j$ in the interval $]x/2,x]$\myp, and relate the
threshold values of $\Delta\nu(\cdot)$ with the lagged ratios
$p_{j+k}/p_j$\myp. This analysis culminates in the proof of
Theorem~\ref{thm:bound}. In Section \ref{sec4}, various sufficient
conditions for $V_n\to\infty$ are derived, which covers the content
of Theorem~\ref{th:->infty}. We also show that these conditions are
not necessary, by constructing examples of weird oscillatory
behavior. In Section \ref{sec5}, we derive a simple integral
condition in terms of the function $\Delta\nu(\cdot)$, necessary and
sufficient in order that $V(t)$ converge to a finite limit. This
criterion is then used to prove Theorem \ref{thm:limit}. In
conclusion, we rehabilitate Karlin's sufficient condition of
convergence, by showing that it is in fact necessary and sufficient.

\section{Poissonization and moment formulas}\label{sec2}
Let $X_{n,j}$ be the occupancy number of box $j$ after $n$ throws,
that is, the number of balls out of the first $n$ that land in
box~$j$. Note that
\begin{equation}\label{Knsum}
K_n=\sum_{j=1}^\infty {\bf 1}\{X_{n,j}>0\}\myp,
\end{equation}
where ${\bf 1}(A)$ is the indicator of event $A$ (i.e., with values
$1$ when $A$ is true and $0$ otherwise). Because $\sum_{j=1}^\infty
X_{n,j}=n$, it is clear that the terms in the sum (\ref{Knsum}) are
not independent.

\subsection{Poissonization}
A common recipe to circumvent the dependence (see \citep{Aldous,
Jacquet-Szpankowski} for a general introduction and \citep{Karlin,
Johnson-Kotz, Kolchin-Sevastyanov-Chistyakov, Hwang-Janson} for
details in the occupancy problem context) is to consider a closely
related model in which the balls are thrown at the jump times of a
unit rate Poisson process $(N(t),\ t\ge 0)$: by this randomization
the balls appear in boxes according to independent Poisson processes
$X_j(t)$, with rate $p_j$ for box $j$. Further advantage of the
poissonized model is that the normalization \eqref{add1} can be
replaced by a weaker summability condition
$\|p\|\equiv\sum_{j=1}^\infty p_j<\infty$, thus allowing one to
avoid computing normalization constants in expressions for $p_j$.
Clearly, the normalization (\ref{add1}) can always be maintained by
rescaling the frequencies $p_j\mapsto \|p\|^{-1}p_j$, to the effect
of a linear time change, $t\mapsto \|p\|\myp t$.

In what follows, we adopt the convention that quantities derived
from the poissonized version of the occupancy problem are written as
functions of the continuous time parameter~$t$, while for the
original model we preserve the notation with lower index $n$. In
particular, we write $X_j(t)$ (cf.\ above) for the number of balls
that land in box $j$ by time $t$ and
\begin{equation}\label{eq:defKt}
K(t):=K_{N(t)} = \sum_{j=1}^\infty {\bf 1}\{X_j(t)>0\}
\end{equation}
for the number of boxes discovered by the Poisson process $N(t)$.
Likewise, denoting by $K_{n,\myp r}$ the number of boxes, each of
which is hit by exactly $r$ of the first $n$ balls, we write
$$
K_r(t):=K_{N(t),\myp r}= \sum_{j=1}^\infty {\bf 1}\{X_j(t)=r\}
$$
for the corresponding poissonized quantity (which is the number of
boxes containing exactly $r$ balls each by time $t$). Clearly,
\begin{equation}\label{eq:KandK}
\begin{aligned}
K_n&=\sum_r K_{n,\myp r}\myp,&\quad K(t)&=\sum_r K_r(t)\myp,\\[.1pc]
n&=\sum_r rK_{n,\myp r}\myp,&\quad N(t)&= \sum_r rK_r(t)\myp.
\end{aligned}
\end{equation}

For the mean values of the number of occupied boxes we have the
formulas
\begin{align}
\label{eq:Phin0}
\Phi_n&:=\EE(K_n)=\sum_{j=1}^\infty (1-(1-p_j)^n)\myp,\\
\label{eq:Phit0}
\Phi(t)&:=\EE(K(t))=\sum_{j=1}^\infty
(1-\e^{-tp_j})\myp,
\end{align}
related by the poissonization identity
$$
\Phi(t)=\e^{-t}\sum_{n=0}^\infty \frac{t^n}{n!}\,\Phi_n\,,
$$
where $\Phi_0=0$. Encoding the collection of frequencies into an
infinite counting measure on $\RR_+\!={}]0,\infty[$
\begin{equation}\label{eq:nu}
\nu(\diff x):=\sum_{j=1}^\infty \delta_{p_j}(\diff x)
\end{equation}
(where $\delta_x$ is the Dirac mass at $x$, i.e., 
$\delta_x(A)={\bf 1}\{x\in A\}$ for $A\subset \RR_+$), we can represent the mean
values (\ref{eq:Phin0}), (\ref{eq:Phit0}) in an integral form as
\begin{align}\label{eq:Phiintn}
\Phi_n&=\int_0^1 \bigl(1-(1-x)^n\bigr)\,\nu(\diff x)\myp,\\
\label{eq:Phiintt} \Phi(t)&=\int_0^\infty (1-\e^{-tx})\,\nu(\diff
x)\myp.
\end{align}
\begin{remark}
When the frequencies are normalized by (\ref{add1}) then all
$p_j\le1$ and the integral in (\ref{eq:Phiintt}) could be written in
the limits from $0$ to $1$, similarly to (\ref{eq:Phiintn}). In the
poissonized model, specific normalization is not important, so we
prefer to use a more flexible notation as in (\ref{eq:Phiintt}). The
same convention applies to similar representations below (see, e.g.,
formulas (\ref{eq:Phitr}) and (\ref{vardiff})).
\end{remark}

Furthermore, set
\begin{align}\label{eq:Phinr}
\Phi_{n,\myp r}&:=\EE(K_{n,\myp r})=
\left(\genfrac{}{}{0pt}{0}{n}{r}\right)\!\int_0^1
x^r(1-x)^{n-r}\,\nu(\diff x),\\[.3pc]
\label{eq:Phitr} \Phi_r(t)&:=\EE[K_{r}(t)]=\frac{t^r}{r!}
\int_0^\infty x^r \e^{-tx} \,\nu(\diff x)\myp,
\end{align}
the latter being related to the derivatives of $\Phi(t)$ via
%\begin{equation}\label{Phider}
$$
\Phi_r(t)=(-1)^{r+1}\, \frac{t^r}{r!}\,\Phi^{(r)}(t)\myp.
$$
%\end{equation}
Note that equations (\ref{eq:KandK}) imply
\begin{equation}\label{eq:PhiandPhi}
\begin{aligned}
\Phi_n&=\sum_r\mypp \Phi_{n,\myp r}\myp,&\quad \Phi(t)&=\sum_r\mypp
\Phi_r(t)\myp,\\[.1pc]
n&=\sum_r r\mypp \Phi_{n,\myp r}\myp,&\quad t&= \sum_r r\mypp
\Phi_r(t)\myp.
\end{aligned}
\end{equation}

An analyst will recognize in (\ref{eq:Phiintt}) a Bernstein function
(see \citep{Berg-Christensen-Ressel}) with the following general
properties (see also \citep{Gnedin-Pitman}).
\begin{lemma}\label{Bernstein}
If an infinite measure\/ $\nu$ on $\RR_+$ satisfies\/
$\int_0^\infty(1-\e^{-x})\,\nu(\diff x)<\infty$\myp\textup{,} then
\textup{(\ref{eq:Phiintt})} defines a function\/ $\Phi(\cdot)$ which
\begin{enumerate}
\item[\rm (i)] is analytic in the right half-plane\myp\textup{,}

\item[\rm (ii)] has alternating derivatives\/
$ (-1)^{r+1}\Phi^{(r)}(t)>0$ \,$(t>0)\myp$\textup{,}

\item[\rm (iii)]
satisfies\/ $\Phi(t)\uparrow\infty$ but\/ $\Phi(t)=o(t)$ as\/
$t\to\infty$\myp.
\end{enumerate}
Conversely\myp\textup{,} if\/ a function\/ $\Phi(t)$ on\/
$[0,\infty[$ has the properties\/ \textup{(ii)} and\/ \textup{(iii)}
along with\/ $\Phi(0)=0$\textup{,} then there exists a unique
infinite measure\/ $\nu$ on\/ $\RR_+$ such that\/ representation
\textup{(\ref{eq:Phiintt})} holds.
\end{lemma}

\subsection{The variance of the number of occupied boxes}\label{sec2.2}
By the independence of summands in (\ref{eq:defKt}), the variance of
$K(t)$ is given by
\begin{equation}\label{eq:V}
V(t):=\Var(K(t))=\sum_{j=1}^\infty (\e^{-tp_j}-\e^{-2tp_j})\myp,
\end{equation}
which is the same as
\begin{equation}\label{vardiff}
V(t)=\int_0^\infty (\e^{-tx}-\e^{-2tx})\,\nu(\diff
x)=\Phi(2t)-\Phi(t)\myp.
\end{equation}

\begin{example}\label{ex:geometric1/2}
For geometric frequencies of ratio $q=1/2$, that is, $p_j=2^{-j}$
($j=1,2,\dots$), the sum (\ref{eq:V}) is evaluated explicitly thanks
to telescoping of partial sums (see \cite[page~1258]{Dutko}):
\begin{align*}
V(t)= \lim_{M\to\infty}\sum_{j=1}^M \left(\e^{-t\myp
2^{-j}}-\e^{-t\myp 2^{-j+1}}\right) = \lim_{M\to\infty}
\bigl(\e^{-t\myp 2^{-M}}-\e^{-t}\myp\bigr) =1-\e^{-t}.
\end{align*}
In particular, it follows that $V(t)\to 1$ as $t\to\infty$. More
generally, a similar simplification occurs in the geometric case
with the ratio $q=2^{-1/k}$ (${k\ge 1}$), where it is convenient to
split the sum in (\ref{eq:V}) into $k$ sub-sums (over
$j=i+k(\ell-1)$, where $i=1,\dots,k$, \,$\ell=1,2,\dots$), each
involving a (non-normalized) geometric sequence with ratio $1/2$.
Applying the previous result (with $q=1/2$) and adding up the $k$
unit contributions emerging in the limit from the $k$ constituent
subsequences, we obtain the convergence $V(t)\to k$ as $t\to\infty$.
For other values of $q$ the formula for the variance does not
simplify.
\end{example}

\begin{remark}
The poissonized variance is additive: if $\bigl(p_j^{(1)}\bigr)$ and
$\bigl(p_j^{(2)}\bigr)$ are two summable sequences of frequencies,
and if $(p_j)$ is obtained by merging them into a single sequence,
then the corresponding variances satisfy
$V^{(1)}(t)+V^{(2)}(t)=V(t)$. This explains the structural
decomposition of the variance mentioned in the Introduction and
illustrated in Example~\ref{ex:geometric1/2}.
\end{remark}

The fixed-$n$ counterpart of (\ref{eq:V}) is
\begin{equation}\label{cross}
V_n =\Phi_{2n}-\Phi_n+ \sum_{i\ne j}^\infty \bigl((1-p_i-p_j)^n
-(1-p_i)^n(1-p_j)^n\bigr)\myp,
\end{equation}
where the cross-terms arise due to dependence in (\ref{Knsum}).

\subsection{Depoissonization}
According to \citep[Proposition 4.3(ii)]{Hwang-Janson}, the
variances $V(n)$ and $V_n$ are always of the same order,
\begin{equation}\label{HJ}
0<\liminf_{n\to\infty}\frac{V(n)}{V_n}\le
\limsup_{n\to\infty}\frac{V(n)}{V_n}<\infty\myp.
\end{equation}

In the next lemma, we establish estimates for the deviation of the
poissonized quantities from their fixed-$n$ counterpart in terms of
higher-order moments, which will be instrumental for
depoissonization in the case of bounded variance (see
Section~\ref{sec3}).
\begin{lemma}\label{varest}
If\/ the normalization\/ \eqref{add1} holds then
\begin{align}
\label{ee1}
\Phi(n)-\Phi_n&=O(n^{-1})\,\Phi_2(n)\myp,\\[.3pc]
\label{ee3}
V(n)-V_n&=O(n^{-1})\,\bigl(\Phi_1(n)^2+\Phi_2(n)\bigr)\myp,
\end{align}
and for each $r=1,2,\dots$
\begin{equation}\label{ee2}
\Phi_r(n)-\Phi_{n,\myp r}=
O(n^{-1})\,\bigl(\Phi_r(n)+\Phi_{r+1}(n)+\Phi_{r+2}(n)\bigr)\myp,
\end{equation}
\end{lemma}

\proof We shall need the elementary inequalities
\begin{equation}\label{eq:<exp<}
0\le \e^{-nx}-(1-x)^n\le n\myp x^2\myp \e^{-nx}\qquad (0\le x\le
1)\myp.
\end{equation}
The first inequality is obvious, while the second one follows from the estimate
$$
(1-x)^n\ge (1-x^2)^n\myp \e^{-nx}\ge (1-n\myp x^2)\mypp
\e^{-nx}\myp.
$$

Now, using representations (\ref{eq:Phiintn}), (\ref{eq:Phiintt})
(rewriting the integral (\ref{eq:Phiintt}) in the limits from $0$ to
$1$, due to \eqref{add1}) and inserting the bounds (\ref{eq:<exp<}),
we obtain
\begin{align*}
0\le \Phi_n-\Phi(n)&=\int_0^1
 \bigl(\e^{-nx}-(1-x)^n\bigr)\,\nu(\diff x)
 \le
\frac{2}{n}\,\Phi_2(n)\myp,
\end{align*}
which proves (\ref{ee1}). Next, from (\ref{eq:Phinr}) and
(\ref{eq:Phitr}) we get
\begin{equation}\label{eq:Phi-Phi}
\Phi_r(n)-\Phi_{n,\myp r} =O(n^{-1})\,\Phi_r(n)+ \frac{n^r}{r!}
\int_0^1 x^r \bigl(\e^{-nx}-(1-x)^{n-r}\bigr)\,\nu(\diff x)\myp.
\end{equation}
By the inequalities (\ref{eq:<exp<}), for each $x\in[0,1]$
\begin{align}\label{eq:lower}
\e^{-nx}-(1-x)^{n-r}&\ge \e^{-nx}-\e^{-(n-r)\myp x}
\ge -(\e^r-1)\myp x\mypp \e^{-nx}\myp,\\[.3pc]
\label{eq:upper}
\e^{-nx}-(1-x)^{n-r}&\le \e^{-nx}-(1-x)^{n}\le
n\myp x^2\myp \e^{-nx}\myp.
\end{align}
Substituting the estimates (\ref{eq:lower}) and (\ref{eq:upper})
into (\ref{eq:Phi-Phi}) and recalling the notation
(\ref{eq:Phitr}) yields~(\ref{ee2}).

Finally, as shown in \citep[Theorem 2.3]{Hwang-Janson}, the
cross-terms in (\ref{cross}) can be evaluated as
\begin{align*}
(1-p_i)^n (1-p_j)^n-(1-p_i-p_j)^n
&= n\myp p_i\myp p_j\mypp (1-p_i)^{n-1}(1-p_j)^{n-1}\\[.3pc]
&\quad+O\bigl(n^2 p_i^2 p_j^2\mypp
(1-p_i)^{n-2}(1-p_j)^{n-2}\bigr)\myp.
\end{align*}
Inserting this estimate into (\ref{cross}) and summing over all
$i,j$, we obtain
\begin{equation}\label{eq:V_n_est}
V_n=\Phi_{2n}-\Phi_n+O(n^{-1})\,\Phi_{n,1}^2+
O(n^{-2})\,\Phi_{n,\myp 2}^2\,.
\end{equation}
From (\ref{eq:PhiandPhi}) and (\ref{eq:Phiintn}) it follows that if
the condition (\ref{add1}) holds then
$$
\Phi_{n,\myp r}\le \Phi_n=\int_0^1 \bigl(1-(1-x)^n\bigr)\,\nu(\diff
x) \le \int_0^1 n\myp x\,\nu(\diff x) =n\myp,
$$
and similarly, using (\ref{eq:Phiintt}),
$$
\Phi_r(n)\le \Phi(n)=\int_0^\infty (1-\e^{-nx})\,\nu(\diff x)\le
\int_0^\infty n\myp x\,\nu(\diff x) =n\myp.
$$
Hence, subtracting (\ref{eq:V_n_est}) from (\ref{vardiff}) and using
the estimates (\ref{ee1}) and (\ref{ee2}), we arrive at~(\ref{ee3}).
\endproof

\section{Bounded variance}\label{sec3}
In this section, we mainly focus on the situation where the variance
$V(t)$ is bounded.

\subsection{Auxiliary estimates}
We first derive various useful inequalities involving the functions
$V(t)$, $\Phi(t)$, $\Phi_r(t)$ and the measure $\nu$\myp. Since
$\Phi'(t)$ is decreasing and $V(t)=\Phi(2t)-\Phi(t)$, the mean value
theorem yields
$$
 \Phi'(2t)\le\frac{\Phi(2t)-\Phi(t)}{t}=\frac{V(t)}{t}\le\Phi'(t)\myp,
$$
or equivalently
\begin{equation}\label{PhiV}
 \frac{1}{2}\,\Phi_1(2t)\le V(t)\le \Phi_1(t)\myp.
\end{equation}
The first inequality in (\ref{PhiV}) generalizes.

\begin{lemma}\label{upbound}
For\/ $r=1,2,\dots$ and\/ $t>0$\textup{,}
$$
\Phi_r(t)\le \frac{2^{r(r+1)/2}}{r!} \,V(2^{-r}t)\myp.
$$
\end{lemma}
\proof Setting $f_r(t):=(-1)^{r+1}\Phi^{(r)}(t)>0$ (see Lemma
\ref{Bernstein}(ii)), we shall prove by induction the equivalent
inequality
\begin{align}\label{eq:Phi_n}
f_r(t)\le \frac{2^{r(r+1)/2}\,V(2^{-r}t)}{t^r} \qquad (t>0)\myp.
\end{align}
Suppose (\ref{eq:Phi_n}) holds for $f_1,\dots,f_{r-1}$. Note that
$f''_{r-1}(t)=f_{r+1}(t)>0$, hence the function $f_{r-1}$ is convex
and therefore
\begin{equation}\label{eq:phi_r1}
\frac{f_{r-1}(t/2)-f_{r-1}(t)}{t/2}\ge -f'_{r-1}(t)=f_r(t)\myp.
\end{equation}
On the other hand, since $f_{r-1}(t)\ge0$ and by the induction
hypothesis,
\begin{equation}\label{eq:phi_r2}
\frac{f_{r-1}(t/2)-f_{r-1}(t)}{t/2} \le\frac{f_{r-1}(t/2)}{t/2} \le
\frac{2^{r(r-1)/2}\,V(2^{-r}t)}
 {(t/2)^r}\mypp.
\end{equation}
Combining (\ref{eq:phi_r1}) and (\ref{eq:phi_r2}), we obtain
(\ref{eq:Phi_n}) for $f_r$. Thus, the induction step follows, and
the proof is complete.
\endproof

Consider the limits superior
\begin{equation}\label{eq:barv}
\bar{v}:=\limsup_{t\to\infty}\myp V(t)\myp,\qquad\bar{\varphi}_r:=
\limsup_{t\to\infty}\mypp\Phi_r(t)\quad(r=1,2,\dots)\myp.
\end{equation}
By continuity, $V(t)$ is uniformly bounded on $[0,\infty[$ if and
only if $\bar{v}<\infty$\myp, and the same is true for $\Phi_r(t)$
in terms of the condition $\bar{\varphi}_r<\infty$.

Note that $\bar{v}$ is strictly positive (cf.\
\citep[page~1258]{Dutko}); indeed, setting $t=1/p_k$ in (\ref{eq:V})
we have
\begin{equation}\label{eq:vbar>0}
\bar{v} \ge \limsup_{k\to\infty} \sum_{j=1}^\infty
\bigl(\e^{-p_j/p_k}-\e^{-2p_j/p_k}\bigr)\ge \e^{-1}-\e^{-2}>0\myp.
\end{equation}

\begin{corollary}\label{cor:v_phibar}
The conditions\/ $\bar{v}<\infty$ and\/ $\bar{\varphi}_1<\infty$ are
equivalent and imply\/ $\bar{\varphi}_r<\infty$ for all\/ $r\ge 1$.
\end{corollary}
\proof Follows from (\ref{PhiV}) and Lemma \ref{upbound}.
\endproof

Appealing to Lemma \ref{varest}, we have depoissonization in terms
of moments.

\begin{corollary}\label{depoissmom}
If\/ $\bar{v}<\infty$ then\myp\textup{,} as\/
$n\to\infty$\myp\textup{,}
$$
\Phi(n)-\Phi_n=O(n^{-1})\myp,\qquad V(n)-V_n=O(n^{-1})\myp,
$$
and\/\textup{,} for all $r\ge1$\/\textup{,}
$$
\Phi_r(n)-\Phi_{n,\myp r}=O(n^{-1})\myp.
$$
\end{corollary}

\subsection{Uniform upper bounds for $\bar{\varphi}_r$}
Lemma \ref{upbound} entails an estimate of $\bar{\varphi}_r$ through
either $\bar{v}$ or $\bar{\varphi}_1$. With some more effort, we
will derive an improved upper bound that does not depend on $r$.
Recall that the measure $\nu$ is defined in (\ref{eq:nu}), and
consider the new (finite) measure
\begin{equation}\label{eq:defnut}
\tilde{\nu}(\diff x):=x\myp \nu(\diff x)=\sum_{j=1}^\infty
p_j\mypp\delta_{p_j}(\diff x)\myp.
\end{equation}
When the normalization (\ref{add1}) holds, this is a probability
measure governing the frequency distribution of the random box
discovered by ball~$1$.

Using the measure $\tilde{\nu}$, we can rewrite (\ref{eq:Phitr}) as
follows
\begin{equation}\label{eq:Phitr-nutilde}
\Phi_r(t)= \frac{t^r}{r!}\int_0^\infty x^{r-1}\mypp \e^{-xt}\,
\tilde{\nu}(\diff x)\myp.
\end{equation}
Also, let us set
\begin{equation}\label{eq:tilde_eta}
\bar{\eta}:=\limsup_{x\downarrow 0}\frac{\tilde{\nu}[0,x]}{x}\mypp.
\end{equation}

\begin{lemma}\label{lm:limsupPhir}
Suppose that\/ $\bar{v}<\infty$. Then for all\/ $r=1,2,\dots$
\begin{equation}\label{eq:phi_vbar}
\bar{\varphi}_r\le \bar{\eta}\le \e\myp \bar{\varphi}_1 \le 2\myp
\e\myp \bar{v}\myp.
\end{equation}
\end{lemma}

\proof Note that the last inequality in (\ref{eq:phi_vbar}) follows
from (\ref{PhiV}). Further, integrating by parts in
(\ref{eq:Phitr-nutilde}) and using the substitution $y=xt$, we get
\begin{equation}\label{eq:Phir}
\Phi_r(t)= \frac{t}{r!}\int_0^\infty \e^{-y}\myp
y^{r-2}\left(y+1-r\right)\tilde{\nu}[0,y/t]\,\diff y\myp.
\end{equation}
For $r=1$, due to monotonicity of the function
$\tilde{\nu}[0,\cdot\mypp]$, (\ref{eq:Phir}) implies
\begin{equation}\label{eq:tilde_rho}
\Phi_1(t)\ge t\int_{1}^\infty \e^{-y}\,\tilde{\nu}[0,y/t]\,\diff y
\ge \e^{-1}\mypp\frac{\tilde{\nu}[0,1/t]}{1/t}\mypp,
\end{equation}
and by letting here $t\to\infty$ we obtain $\bar{\varphi}_1\ge
\e^{-1}\bar{\eta}\myp$ (see
(\ref{eq:tilde_eta}),~(\ref{eq:phi_vbar})).

On the other hand, for any $r\ge 1$ from (\ref{eq:Phir}) it follows
that
$$
\Phi_r(t)\le \frac{1}{r!}\int_0^\infty \e^{-y}\myp
y^{r}\mypp\frac{\tilde{\nu}[0,y/t]}{y/t}\,\diff y\qquad (r\ge
1)\myp,
$$
which implies \myp$\bar{\varphi}_r\le \bar{\eta}$\myp{} by the
``$\limsup$" part of Fatou's lemma \citep[\S\mypp
IV.2\myp]{Feller2}.
\endproof

%\textit{Health \& Safety}: Let us justify the use of Fatou's lemma.
%Note that the function $\eta(u):=u^{-1}\tilde{\nu}[0,u]$ is bounded on
%\RR_+$, since $\lim_{u\to\infty}\eta(u)=0$ (because $\tilde{\nu}[0,u]\le
%\|p\|=\sum_j p_j<\infty$), while
%$\limsup_{u\downarrow0}\eta(u)=\bar{\eta}<\infty$. Hence, the
%function under the integral is uniformly bounded by the (integrable)
%function $cy^{r}\e^{-y}$, and we can apply the $\limsup$ Fatou's
%lemma:
%$$
%\bar{\varphi}_r=\limsup_{t\to\infty}\Phi_r(t)\le
%\frac{1}{r!}\int_0^\infty y^{r}\e^{-y}\limsup_{t\to\infty}
%\eta(y/t)\,\diff u=\frac{\bar{\eta}}{r!\ref{eq:(i)})}\int_0^\infty
%y^{r}\e^{-y}\,\diff u=\bar{\eta}\myp.
%$$

\subsection{Growth of the mean number of occupied boxes}
Lemma \ref{lm:limsupPhir} implies that if $\bar{v}<\infty$ then each
term in the decomposition $\Phi(t)= \sum_{r=1}^\infty \Phi_r(t)$
makes a uniformly bounded contribution to $\Phi(t)\to\infty$. This
is to be contrasted with the case of frequencies akin to $p_j\sim
c\myp j^{-\alpha}$ ($\alpha>1$), where $V(t)$, $\Phi(t)$ and
$\Phi_r(t)$ ($r\ge1$) are of the same order $O(t^\alpha)$ as
$t\to\infty$ (see \citep{Karlin}). The next lemma estimates the
growth of $\Phi(t)$ in the case of bounded variance.

\begin{lemma}\label{lm:Phi<V}
Suppose that\/ $\bar{v}<\infty$\myp. Then
$$
 \limsup_{t\to\infty}\frac{\Phi(t)}{\log t}\le 2\bar{v}\myp.
$$
\end{lemma}
\proof For any $\varepsilon>0$\myp, there exists $t_0>0$ such that
for all $t\ge t_0$
$$
\Phi_1(t)\le \bar\varphi_1+\varepsilon\le 2\bar{v} +\varepsilon\myp,
$$
due to Lemma \ref{lm:limsupPhir}. Therefore,
$$
\Phi(t)-\Phi(t_0)=\int_{t_0}^t\Phi'(s)\,\diff s
=\int_{t_0}^t\frac{\Phi_1(s)}{s}\,\diff s\le (2\bar{v}
+\varepsilon)\myp(\log t-\log t_0)\myp.
$$
Hence,
$$
\limsup_{t\to\infty}\frac{\Phi(t)}{\log t}=
\limsup_{t\to\infty}\frac{\Phi(t)-\Phi(t_0)}{\log t-\log t_0}\le
 2\bar{v} +\varepsilon\myp,
$$
and since $\varepsilon>0$ is arbitrary, our claim follows.

A shorter proof is by a simple ``$\limsup$'' version of
L'H\^opital's rule:
%L'Hospital's
$$
\limsup_{t\to\infty}\frac{\Phi(t)}{\log t}\le
\limsup_{t\to\infty}\frac{\Phi'(t)}{1/t}=\limsup_{t\to\infty}
\Phi_1(t)=\bar\varphi_1\le 2\myp\bar v\myp,
$$
due to Lemma \ref{lm:limsupPhir}.
\endproof

\subsection{The basic representation of the variance $V(t)$}
As in \citep{Karlin}, it is convenient to rewrite the formula
(\ref{vardiff}) for the variance  as a single integral
representation. Recall that $\nu$ is given by (\ref{eq:nu}), and
introduce the function
\begin{equation}\label{eq:Delta_nu}
\Delta\nu(x):=\nu\mypp]x/2,\,x]=\#\{j: x/2<p_j\le
x\}\qquad(x>0)\myp.
\end{equation}

\begin{lemma}\label{lm:V}
The variance\/ $V(t)$ can be represented as
\begin{equation}\label{eq:V1}
V(t)=t\int_0^\infty \e^{-tx} \Delta\nu(x)\,\diff x\qquad(t\ge0)\myp.
\end{equation}
\end{lemma}
\proof Setting $\nu_{\rm c}(x):=\nu\mypp]x,\infty[$ and integrating
by parts in (\ref{vardiff}) gives
\begin{align*}
V(t)&=\int_0^\infty (\e^{-2tx}-\e^{-tx})\,\diff\nu_{\rm c}(x)\\
&=\bigl(\e^{-2tx}-\e^{-tx}\bigr)\mypp\nu_{\rm c}(x)\bigr|_0^\infty
+t\int_0^\infty \e^{-tx}\bigl(\nu_{\rm c}(x/2)-\nu_{\rm c}(x)\bigr)\,\diff x\\
&=t\lim_{x\downarrow 0} x\myp\nu_{\rm c}(x) +t\int_0^\infty
\e^{-tx}\Delta\nu(x)\,\diff x\myp,
\end{align*}
and (\ref{eq:V1}) will follow if we show that $x\myp\nu_{\rm
c}(x)\to0$ as $x\downarrow0$. To this end, note that the mean value
of the measure $\nu$ is finite: $\int_0^\infty x\,\nu(\diff
x)=\sum_{j=1}^\infty p_j<\infty$\myp. Hence, integration by parts
yields
\begin{equation}\label{eq:by_parts}
\infty>\int_0^\infty x\mypp\nu(\diff x) =\lim_{x\downarrow 0}
x\myp\nu_{\rm c}(x)+\int_0^\infty \nu_{\rm c}(x)\,\diff x\myp,
\end{equation}
and it follows that the limit in (\ref{eq:by_parts}) exists and,
moreover, must vanish, for otherwise
%the function $\nu_{\rm c}(x)$ would behave like $c\mypp x^{-1}$
%for small $x$, in which case
the integral on the right-hand side of (\ref{eq:by_parts}) would
diverge.
\endproof

\begin{corollary}\label{cor:D}
The function
\begin{equation}\label{eq:D}
D(x):=\int_0^x \Delta\nu(u)\,\diff u
\end{equation}
is well defined and uniformly bounded for all $x\ge0$\/. In
particular\myp\textup{,} $D(0)=0$\/.
\end{corollary}
\proof Letting $t=1$ in (\ref{eq:V1}), we obtain
$$
V(1)\ge \int_0^x \e^{-u} \Delta\nu(u)\,\diff u \ge \e^{-x} \int_0^x
\Delta\nu(u)\,\diff u\myp,
$$
hence $D(x)\le \e^{x}\myp V(1)<\infty$ for any $x>0$\myp. Vanishing
at zero is obtained by the absolute continuity of the integral.
Finally, boundedness of $D(x)$ follows because $\Delta\nu(x)\equiv
0$ for all $x$ large enough.
\endproof

Integrating by parts in (\ref{eq:V1}) and using Corollary
\ref{cor:D}, we obtain an alternative representation, which will
also be useful:
\begin{equation}\label{eq:V-D}
V(t)=t^2\int_0^\infty \e^{-tx} D(x)\,\diff x=\int_0^\infty
\e^{-y}\myp y \,\frac{D(y/t)}{y/t}\,\diff x\qquad(t>0)\myp.
\end{equation}

\subsection{Estimates using the function $\Delta\nu(x)$}\label{sec3.5}

It is immediately clear from (\ref{eq:V1}) that if $\Delta\nu(x)\le
c$ for all $x>0$ then $V(t)\le c$ for all $t>0$. Moreover, one can
obtain two-sided asymptotic bounds as follows.
\begin{lemma}\label{lm:bound}
Recall that\/ $\bar v$ is given by\/
\textup{(\ref{eq:barv})}\textup{,} and set
%\begin{equation}\label{eq:bard}
$$
\bar{w}:=\limsup_{x\downarrow 0}\mypp \Delta\nu(x)\myp.
$$
%\end{equation}
Then\/ $\bar{v}<\infty$ if and only if\/ $\bar{w}<\infty$\textup{,}
and in this case
\begin{equation}\label{eq:v<d<v}
 (\sqrt{5}-2)\,\bar{w}\le \bar{v}\le\bar{w}\mypp.
\end{equation}
\end{lemma}
\proof The substitution $y=tx$ in (\ref{eq:V1}) yields
%%\begin{equation}\label{eq:V2}
$$
V(t)=\int_0^\infty \e^{-y}\Delta\nu(y/t)\,\diff y\myp,
$$
%\end{equation}
and an application of the ``$\limsup$'' part of Fatou's lemma
\citep[\S\mypp{}IV.2]{Feller2} implies
$$
\bar v\le \bar{w} \int_0^\infty \e^{-y}\,\diff y=\bar{w}\mypp.
$$
%\textit{Health \& Safety}: The use of Fatou's lemma in (\ref{eq:V2})
%is justified by the same argument as before, this time using that
%$\Delta\nu(\cdot)$ is finitely supported (more precisely,
%$\Delta\nu(u)=0$ for all $u\ge 2\|p\|$).

For the converse inequality, we need to exploit the special
structure of the measure~$\nu$. Fixing $x>0$ and retaining in
(\ref{eq:V}) the terms with $p_j\in{}]x/2,x]$ only, we obtain
\begin{equation}\label{eq:Vn1}
V(t)\ge \Delta\nu(x) \min _{p\in[x/2,\myp x]}
\bigl(\e^{-tp}-\e^{-2tp}\bigr).
\end{equation}
It is clear that the minimum in (\ref{eq:Vn1}) is attained at one of
the endpoints, that is, $p=x/2$ or $p=x$. Setting
$y=\e^{-tx/2}\in[0,1]$, we note that
$$
\min{}\{y-y^2, y^2-y^4\}=\left\{\begin{array}{ll}
y^2-y^4,&0\le y\le \phi\myp,\\[.2pc]
y-y^2,&\phi\le y\le 1\myp,
\end{array}\right.
$$
where $\phi=(\sqrt{5}-1)/2$ is the golden ratio, which appears here
as the root of the equation $y^2-y^4=y-y^2$ on $]0,1[$\mypp. It is
then easy to see that the right-hand side of (\ref{eq:Vn1}), as a
function of $t$, attains its maximum value $\phi-\phi^2=\sqrt{5}-2$
at $t(x)=2\myp x^{-1}\log{}(1/\phi)\to\infty$ \,($x\downarrow 0$).
Hence $V(t(x))\ge (\sqrt{5}-2)\myp\Delta\nu(x)$, and the first
inequality in (\ref{eq:v<d<v}) follows.
\endproof

Our next goal is to characterize the link between the upper (lower)
bounds on the values of the function $\Delta\nu(x)$ (for small~$x$)
and the lagged frequency ratios $p_{j+k}/p_j$ (for large $j$) with
regard to the threshold value $1/2$.

\begin{lemma}\label{p-and-nu}
For a given positive integer $k$\/\textup{,} the bound
\begin{equation}\label{Deltanuk}
\Delta\nu(x)\le k
\end{equation}
is valid for all sufficiently small\/ $x>0$ if and only if the
condition
\begin{equation}\label{eq:boundrat1}
\frac{p_{j+k}}{p_j}\le \frac{1}{2}
\end{equation}
is satisfied for all sufficiently large\/ $j$\/. The similar
assertion holds true when the sign\/ $\le$ in both\/
\textup{(\ref{Deltanuk})} and\/ \textup{(\ref{eq:boundrat1})} is
replaced by\ $\ge$\mypp.
\end{lemma}

\proof The first part of the lemma (i.e., with \,$\le$\mypp) is just
a reformulation of definitions (see~(\ref{eq:Delta_nu})). Indeed,
applying (\ref{Deltanuk}) with $x=p_j$ implies $p_{j+k}\le
p_j/2$\myp, which is (\ref{eq:boundrat1}). Conversely, if $p_j\le x<
p_{j-1}$ then by (\ref{eq:boundrat1}) we have $p_{j+k}\le p_j/2\le
x/2$, and hence $\Delta\nu(x)=\nu\mypp]x/2,x]\le k$ as required by
(\ref{Deltanuk}).

The ``mirror'' part (i.e., with \,$\ge$\mypp) needs a bit more care.
First, note that it suffices to prove the ``only if'' statement in
the case where $p_j>p_{j+1}$, for if $p_j=p_r$ ($r>j$) then
$p_{j+k}/p_j \ge p_{r+k}/p_{r}$\myp. Now, if $x\in[p_{j+1},p_j[$
then the condition $\Delta\nu(x)\ge k$ implies that
$p_{j+k}>x/2$\myp, whence by letting $x\uparrow p_j$ we get
$p_{j+k}\ge p_j/2$\myp. Similarly, the ``if'' part follows by noting
that $p_{j+k}\ge p_j/2$ implies $\Delta\nu(x)\ge k$ for each
$x\in[p_{j+1},p_j[$\myp.
\endproof

\subsection{Refined asymptotic estimates}\label{sec3.6}
By Lemma \ref{p-and-nu} and the inequality (\ref{eq:v<d<v}),
the upper bound (\ref{Deltanuk}) implies $\bar{v}\le \bar{w}\le k$.
In some cases, however, such an estimate may not be sharp, as the
next example demonstrates.

\begin{example}\label{ex:j/2}
Let $p_j=j\mypp 2^{-j}\in \RT_{1/2}$, so by Theorem \ref{thm:limit}
we have $\lim_{t\to\infty} V(t)=1$. On the other hand,
(\ref{eq:boundrat1}) holds starting from $k=2$, which leads to the
crude bound $\bar{v}\le 2$. An inspection shows that
$\Delta\nu(\cdot)=1$ on $[2p_{i+1},p_{i-1}[$\ and
$\Delta\nu(\cdot)=2$ on $[p_{i},2p_{i+1}[$ \,($i\ge 4$). For a given
$x\in [p_{j},p_{j-1}[$\,, ``excess'' over the value $1$ on the
interval \myp$]0,x]$ occurs on a set of total Lebesgue's measure
bounded by $\sum_{i\ge j} (2p_{i+1}-p_{i})=\sum_{i\ge j}
2^{-i}=2^{-j+1}$, which is small as compared to $x\ge p_{j}$
\,($j\to\infty$).
\end{example}

This example suggests the following refinement of Lemma
\ref{p-and-nu}.
\begin{lemma}\label{p-and-D}
If for some\/ $k\in\NN$ the frequencies\/ $(p_j)$ satisfy
\begin{equation}\label{eq:boundrat11}
\limsup_{j\to\infty}\frac{p_{j+k}}{p_j}\le \frac{1}{2}\mypp,
\end{equation}
then $\limsup_{t\to\infty} V(t)\le k$\myp. The assertion remains
valid when the symbols\/ \,$\le$\, and\/ \,$\limsup$\, are
simultaneously replaced by\/ \,$\ge$\, and\/ \,$\liminf$.
\end{lemma}

\proof It suffices to assume that $k=1$, as the general case would
then follow by the additivity argument (see the remark after
Example~\ref{ex:geometric1/2}). According to (\ref{eq:boundrat11})
(with $k=1$), for any $\varepsilon\in\myp]0,1/5]$ and all
sufficiently large $i$ we have $p_{i+1}/p_i\le 1/2+\varepsilon$\myp.
Hence, $p_{i+1}/p_{i-1} \le (1/2 + \varepsilon)^2\le 49/100<
1/2$\myp, and Lemma \ref{p-and-nu} implies that $\Delta\nu(x)\le 2$
for all sufficiently small~$x$.

On the other hand, using the definition of the function
$\Delta\nu(\cdot)$ one can check that $\Delta\nu(x)\le 1$ when
$x\in[\myp p_{i}\wedge (2p_{i+1}),p_{i-1}[$\mypp. That is to say,
the value $\Delta\nu(x)=2$ may only occur on a subset of $[\myp
p_{i},p_{i-1}[$ with Lebesgue's measure not exceeding
$(2p_{i+1}-p_{i})\vee 0$ (here $a\wedge b:=\min\{a,b\}$, \,$a\vee
b:=\max\{a,b\}$). Therefore, for $x\in[\myp p_{i},p_{i-1}[$\ we have
$$
\int_{p_i}^{x} \Delta\nu(u)\,\diff u\le x-p_i+ (2p_{i+1}-p_i)\vee
0\le x-p_i+2\varepsilon\myp p_i\myp.
$$
Inserting these estimates into (\ref{eq:D}), we obtain for
$x\in[p_{j},p_{j-1}[$
\begin{align*}
D(x)&=\int_{p_j}^{x} \Delta\nu(u)\,\diff u+\sum_{i>
j}\int_{p_i}^{p_{i-1}} \Delta\nu(u)\,\diff u\\
&\le x+2\varepsilon\sum_{i\ge j}p_i\le  x + 2\varepsilon
p_{j}\sum_{\ell=0}^\infty \left(\frac{1}{2}+\varepsilon\right)^\ell=
x + \frac{4\varepsilon p_j}{1-2\varepsilon}\mypp.
\end{align*}
It follows that
$$
\frac{D(x)}{x}\le 1 + \frac{4\varepsilon p_j}{(1-2\varepsilon)\myp
x}\le 1+\frac{4\varepsilon}{1-2\varepsilon}\to
1\qquad(\varepsilon\to0)\myp,
$$
hence $\limsup_{x\downarrow 0} D(x)/x \le 1$. Finally, applying to
(\ref{eq:V-D}) the ``$\limsup$'' part of Fatou's lemma
\citep[\S\mypp IV.2\myp]{Feller2}, we obtain
%\begin{equation}\label{eq:limsup}
$$
\limsup_{t\to\infty} \int_0^\infty \e^{-y}\mypp
y\,\frac{D(y/t)}{y/t}\,\diff y\le \int_0^\infty \e^{-y}\mypp
y\,\diff y=1\myp,
$$
%\end{equation}
and the first half of the lemma is proved.

For the second half (with $\ge$ and $\liminf$\myp), suppose again
that $k=1$. According to (\ref{eq:boundrat11}), for any
$\varepsilon\in\myp]0,1/2]$ and all sufficiently large $i$, we have
$p_{i}/p_{i-1}\ge 1/2 - \varepsilon$\myp. Observe that possible
deviations of the function $\Delta\nu(\cdot)$ from value $1$ may
only occur as follows: if $p_{i-1}<2p_{i}$ then $\Delta\nu(x)\ge 2$
for $x\in[p_{i-1},2p_{i}[$\myp, while if $p_{i-1}>2p_{i}$ then
$\Delta\nu(x)= 0$ for $x\in[2p_{i},p_{i-1}[$\myp. In either case,
the contribution of the interval with the endpoints $p_{i-1}$ and
$2p_{i}$ to the integral $\int_0^x( \Delta\nu(u)-1)\,\diff u=D(x)-x$
is bounded from below by
$$ \int_{p_{i-1}\wedge (2p_{i})}^{p_{i-1}\vee
(2p_{i})}\left(\Delta\nu(u)-1\right)\diff u\ge
2p_{i}-p_{i-1}\myp.
$$
Using this remark, for a given $x\in {}]p_{j},p_{j-1}]$ we obtain
\begin{align*}
\int_0^x \left(\Delta\nu(u)-1\right)\diff u &\ge
(2p_{j}-p_{j-1})\wedge 0+\sum_{i> j} (2p_{i}-p_{i-1})\\
&=(2p_{j}-p_{j-1})\wedge 0 -p_j+\sum_{i> j} p_{i}\\
&\ge 2p_j\left(1-\frac{p_{j-1}}{2p_j}\right)\wedge 0
-p_j+p_j\sum_{\ell=1}^\infty
\left(\frac{1}{2}-\varepsilon\right)^\ell\\
&\ge \frac{-4p_j\myp\varepsilon}{1-2\varepsilon}\wedge 0
-p_j+\frac{p_j\myp(1-2\varepsilon)}{1+2\varepsilon}\\
&= -\frac{4p_j\myp\varepsilon}{1-2\varepsilon}
-\frac{4p_j\myp\varepsilon}{1+2\varepsilon}=
-\frac{8p_j\myp\varepsilon}{1-4\varepsilon^2}\mypp.
\end{align*}
Hence
$$
\frac{D(x)}{x}\ge 1 - \frac{8p_j\myp\varepsilon
}{(1-4\varepsilon^2)\myp x}\ge
1-\frac{8\varepsilon}{1-4\varepsilon^2}\mypp,
$$
and since $\varepsilon $ is arbitrary, it follows that
$\liminf_{x\downarrow 0} D(x)/x \ge 1$. It remains to use Fatou's
lemma in (\ref{eq:V-D}) to conclude that $\liminf_{t\to\infty}
V(t)\ge 1$\myp.
\endproof

\begin{corollary}\label{cor:lagged}
Suppose that the condition\/ \textup{(\ref{pv})} is satisfied for
some\/ $k\in\NN$\myp\textup{,} that is\myp\textup{,} $p_{j+k}/p_j\to
1/2$ as\/ $j\to\infty$. Then\/ $V(t)\to k$ as\/ $t\to\infty$.
\end{corollary}
\proof Readily follows by combining the two halves of Lemma
\ref{p-and-D}.
\endproof
Note that Corollary \ref{cor:lagged} is exactly the ``if'' part of
Theorem \ref{thm:limit}. In Section~\ref{sec5} below, where the
issue of converging variance is considered in detail, we will give a
direct, shorter proof of the necessity of the condition~(\ref{pv}).
\begin{example}
Note that a converse statement to either half of Lemma \ref{p-and-D}
is \emph{not valid}\/. Indeed, if $(p_j)\in\RT_q$ with
$q\in[0,1/2[$\myp, then $\Delta\nu(\cdot)=1$ on
$[p_{i},2p_{i}[$\myp\ and $\Delta\nu(\cdot)=0$ on
$[2p_{i},p_{i-1}[$\myp\ (for $i$ large enough). This implies that
the graph $y=D(x)/x$ consists of arcs of hyperbolas with alternating
monotonicity (supported on intervals of the form
$[p_{i},2p_{i}[$\myp\ and $[2p_{i},p_{i-1}[$\myp), and in particular
\begin{equation}\label{eq:max/min}
\begin{aligned}
\max_{x\in[p_j,\myp p_{j-1}]}\frac{D(x)}{x}&=\frac{D(2p_j)}{2p_j}
=\frac{1}{2p_j}\sum_{i\ge j} p_i=\frac{1+\rho_j}{2}\myp,\\
\min_{x\in[2p_{j},\myp
2p_{j-1}]}D(x)&=\frac{D(p_{j-1})}{p_{j-1}}=\frac{1}{p_{j-1}}\sum_{i\ge
j} p_i=q\myp(1+\rho_j)\myp,
\end{aligned}
\end{equation}
where $\rho_j=p_j^{-1}\sum_{i>j} p_{i}$ (cf.~(\ref{eq:rho})). The
$\RT_q$\myp-condition implies that $\rho_j\to q/(1-q)$ as
$j\to\infty$, so from (\ref{eq:max/min}) we get
\begin{equation}\label{eq:DD}
\frac{q}{1-q}\le \liminf_{x\downarrow 0} \frac{D(x)}{x}\le
\limsup_{x\downarrow0} \frac{D(x)}{x}\le \frac{1}{2\myp(1-q)}\mypp.
\end{equation}

In particular, setting $q=0$ (e.g., when $(p_j)$ is a Poisson
distribution) and taking a ``doubled'' sequence (i.e., determined by
$\nu(\diff x) = \sum_{j=1}^\infty 2\mypp\delta_{p_j}(\diff x)$)\myp,
by the additivity argument we get $\limsup_{t\to\infty} V(t)\le
2\cdot(1/2)=1$, while $\limsup_{j\to\infty} p_{j+1}/p_j=1$.
Likewise, choosing $q=1/3$ and again doubling the sequence, from
(\ref{eq:DD}) and by Fatou's lemma applied to (\ref{eq:V-D}), we
obtain that $\liminf_{t\to\infty} V(t)\ge 2\cdot(1/2)=1$, whereas
$\liminf_{j\to\infty} p_{j+1}/p_j=q=1/3$\myp.
\end{example}

\subsection{Proof of Theorem \textup{\ref{thm:bound}}}
We are now in a position to prove Theorem \ref{thm:bound}, and let
us start by proving its poissonized version. By Lemma
\ref{lm:bound}, the conditions $\bar{v}<\infty$ and $\bar{w}<\infty$
are equivalent, and according to the first half of Lemma
\ref{p-and-nu}, the latter condition holds if and only if
(\ref{eq:boundrat1}) is satisfied for some $k\in\NN$, which is
equivalent to (\ref{eq:boundsup}) (possibly, with a bigger~$k$).

The second part of the theorem (leading to the estimate $\bar{v}\le
k$) is settled by Lemma~\ref{p-and-D}, since condition
(\ref{eq:boundrat11}) of the lemma coincides with condition
(\ref{eq:boundsup}) of the theorem.

Furthermore, by (\ref{HJ}) the condition $\limsup_{n\to\infty}
V_n<\infty$ is equivalent to $\bar{v}<\infty$, in which case also
$\limsup_{n\to\infty} V_n=\bar{v}$ by Corollary \ref{depoissmom}.

Finally, the optimality of the bound $\bar{v}\le k$ follows by
merging $k$ geometric sequences with ratio $q=1/2$ each and using
the additivity argument (alternatively, one can consider the
geometric frequencies with ratio $q=2^{-1/k}$).

Thus, the proof of Theorem \ref{thm:bound} is complete.

\subsection{Comment on the threshold constant}
Let us remark that the threshold $1/2$ in Theorem \ref{thm:bound} is
chosen to match neatly with Theorem \ref{thm:limit}. Replacing $1/2$
in (\ref{eq:boundsup}) by some other value $0<q<1$ would lead to a
more sophisticated upper bound
\begin{equation}\label{eq:C_k,q}
\limsup_{n\to\infty} V_n \le k\myp \lceil \log_{1/q}2\rceil\myp,
\end{equation}
where $\lceil x\rceil:=\min{}\{m\in\ZZ: m\ge x\}$ is the ceiling
integer part of $x$. Indeed, iterating the condition
$\limsup_{j\to\infty} p_{j+k}/p_j\le q$, we get
$\limsup_{j\to\infty}p_{j+ik}/p_j\le q^i\le 1/2$\myp, provided that $i\ge
\lceil \log_{1/q} 2\rceil$, and (\ref{eq:C_k,q}) follows by
Lemma~\ref{p-and-D}.

In fact, the constant $\lceil \log_{1/q}2\rceil$ here has the
meaning of an upper bound for $\limsup_{n\to\infty} V_n$ in the
geometric case with ratio $q$\myp. Note that the representation
(\ref{Arch}) leads to a similar (in general, slightly better)
estimate $\limsup_{n\to\infty} V_n\le k\mypp\bigl(\log_{1/q} 2+\max
\delta_V(\cdot)\bigr)$ \,(cf.~(\ref{eq:C_k,q})).

\section{Convergence to infinity}\label{sec4}

In this section, we establish new sufficient conditions in order
that $V(t)\to\infty$ as $t\to\infty$ (which, in view of (\ref{HJ}),
is equivalent to $V_n\to\infty$ as $n\to\infty$\myp). Note that the
combination of Theorems \ref{suff1} and \ref{sctail} (to be proved
in Sections \ref{sec4.1} and \ref{sec4.2}, respectively) along with 
discussion in Section~\ref{sec4.4} will settle
Theorem \ref{th:->infty} stated in the Introduction.

\subsection{First set of sufficient conditions}\label{sec4.1} It is natural to
seek a condition for $V(t)\to\infty$ based on the representation
(\ref{eq:V}), that is, in terms of the function $\Delta\nu(x)$. In
turn, such a condition may be transformed into the information about
the lagged ratio $p_{j+k}/p_j$ (cf.\ Theorem \ref{thm:bound}).
\begin{theorem}\label{suff1}
The condition
\begin{equation}\label{eq:(i)}
\smash[b]{\lim_{x\downarrow 0}\Delta\nu(x)=\infty}
\end{equation}
implies that
\begin{equation}\label{eq:(ii)}
\smash[t]{\forall\myp k\in\NN\myp,\ \ \quad \liminf_{j\to\infty}
\frac{p_{j+k}}{p_j}\ge \frac{1}{2}\mypp,}
\end{equation}
which in turn implies that\/ $V(t)\to\infty$ as\/ $t\to\infty$\myp.
\end{theorem}

\proof
If condition (\ref{eq:(i)}) holds then for any $k\in\NN$ we
have $\Delta\nu(x)\ge k$ for all sufficiently small $x>0$. By Lemma
\ref{p-and-nu}, this implies that $p_{j+k}/p_j\ge 1/2$ for all $j$
large enough, and (\ref{eq:(ii)}) follows. Further, condition
(\ref{eq:(ii)}) implies convergence of $V(t)$ to infinity by
Lemma~\ref{p-and-D}\myp.
\endproof

Note that condition (\ref{eq:(ii)}) is obviously fulfilled for any
sequence $(p_j)$ from $\RT_1$, in which case it is well known that
$V(t)\to\infty$ \cite{Karlin, Dutko}. The next example demonstrates
that there are instances of frequencies $(p_j)$ satisfying
(\ref{eq:(i)}) but \emph{not} in $\RT_1$. This example will also
show that conditions (\ref{eq:(i)}) and (\ref{eq:(ii)}) of Theorem
\ref{suff1} are \textit{not necessary} in order that
$V(t)\to\infty$\myp.
\begin{example}\label{ex:3^k}
Let $0<q<1$ and suppose that the sequence $(p_j)$ consists of the
values $q^{i}$, each repeated $i$ times ($i=1,2,\dots$), which
corresponds to the measure $ \nu(\diff x)=\sum_{i=1}^\infty
i\myp\delta_{q^{i}}(\diff x)$\myp. Note that the sequence $(p_j)$ is
not in any $\RT$-class, since $\limsup_{j\to\infty}p_{j+1}/p_j=1$
but $\liminf_{j\to\infty}p_{j+1}/p_j=q$\myp. However, for any
$q\in{}]0,1[$ we have $V(t)\to\infty$\myp, since for
$t\in[q^{-j},\myp q^{-j-1}]$
\begin{align*}
V(t)&=\sum_{i=1}^\infty i\mypp\bigl(\e^{-q^{i}t}-\e^{-2\myp
q^{i}t}\bigr)\ge
j\mypp\bigl(\e^{-q^{j}t}-\e^{-2\myp q^{j}t}\bigr)\\
&\ge j \min_{y\in [1,\mypp
q^{-1}]}\bigl(\e^{-y}-\e^{-2y}\bigr)=j\mypp\bigl(\e^{-1/q}-\e^{-2/q}\bigr)
\to \infty\qquad(j\to\infty)\myp.
\end{align*}
If $1/2\le q<1$ then for $x\in [q^j,q^{j-1}[$\ we have
$\Delta\nu(x)\ge j\to\infty$ as $x\downarrow 0$, and condition
(\ref{eq:(i)}) is valid. On the other hand, if $0<q<1/2$ then
$\Delta\nu(x)=0$ for $x\in [2q^j,q^{j-1}[$\mypp, hence
$\liminf_{x\downarrow 0} \Delta\nu(x)=0$ and (\ref{eq:(i)}) fails.
Also, for any $k\ge 1$, we have
$\liminf_{j\to\infty}p_{j+k}/p_j=q<1/2$\myp, so condition
(\ref{eq:(ii)}) is not valid.
\end{example}

\subsection{Another set of conditions}\label{sec4.2}
A different sufficient condition exploits the link between $V(t)$
and the mean number of singleton boxes $\Phi_1(t)$, as in Lemma
\ref{lm:limsupPhir}. An equivalent condition may be set in terms of
the tail ratio $\rho_j=p_j^{-1}\sum_{i>j}^\infty p_{i}$
(see~(\ref{eq:rho})). Recall the definition (\ref{eq:defnut}) of the
measure $\tilde{\nu}$\myp.

\begin{theorem}\label{sctail}
The condition
\begin{equation}\label{eq:nutoverx}
\smash[b]{\lim_{x\downarrow 0} \frac{\tilde{\nu}[0,x]}{x}=\infty}
\end{equation}
is equivalent to
\begin{equation}\label{eq:pjovertail}
\smash[t]{\lim_{j\to\infty} \rho_j=\infty\myp,}
\end{equation}
and each one implies that\/ $V(t)\to\infty$ as\/ $t\to\infty$\myp.
\end{theorem}
\proof By the estimate (\ref{eq:tilde_rho}), condition
(\ref{eq:nutoverx}) implies $\Phi_1(t)\to\infty$, which is
equivalent to $V(t)\to\infty$ by~\eqref{PhiV}. So it remains to show
that (\ref{eq:nutoverx}) and (\ref{eq:pjovertail}) are equivalent to
each other. Observe that for $p_{j+1}\le x< p_j$ we have
$x^{-1}\tilde{\nu}[0,x]\ge\rho_j$, hence (\ref{eq:pjovertail})
implies (\ref{eq:nutoverx}). To prove the converse, note that if
$p_{j+1}=p_j$ then $\rho_j=1+\rho_{j+1}$, so it suffices to consider
the case where $p_{j+1}<p_j$. Then
$$
\rho_j=\inf_{p_{j+1}\le
x<p_j}\frac{\tilde{\nu}[0,x]}{x}\to\infty\qquad(j\to\infty)\myp,
$$
when the condition (\ref{eq:nutoverx}) holds, and hence
(\ref{eq:pjovertail}) follows.
\endproof

\subsection{A counterexample to Theorem\/ \textup{\ref{sctail}}\/}\label{sec4.3}
We construct here an example demonstrating that conditions
(\ref{eq:nutoverx}), (\ref{eq:pjovertail}) are \textit{not
necessary} in order that $V(t)\to\infty$ (or, equivalently,
$\Phi_1(t)\to\infty$). In particular, due to the estimate
(\ref{eq:phi_vbar}) (with $r=2$), this example will show that
$V(t)\to\infty$ does not necessarily imply $\Phi_2(t)\to\infty$. On
the other hand, in view of the inequality
\begin{align*}
2\myp\Phi_2(t)\ge \sum_{1/2<tp_j\le 1}(tp_j)^2\myp \e^{-tp_j} &\ge
\frac{\e^{-1}}{4}\, \#\bigl\{p_j\in{}]1/(2t),1/t]\bigr\}
=\frac{\e^{-1}}{4}\,\Delta\nu (1/t)\myp,
\end{align*}
it is a priori clear that $\Phi_2(t)$ cannot be uniformly bounded in
such a situation, because $\bar{w}=\infty$ according
to~(\ref{eq:v<d<v}).

\begin{example}\label{ex:Phi1notPhi2}
Let $k_0,k_1,k_2,\dots$ be an increasing integer sequence. Take 
%Choose
the frequencies $(p_j)$ in the form
\begin{equation}\label{eq:p_ex}
p_j=\left\{
\begin{array}{ll}
\displaystyle
k_{1}^{-1},&0<j\le k_0\mypp,\\[.3pc]
\displaystyle k_{i+1}^{-1},\ \ &k_0+\dots+k_{i-1}<j\le
k_0+\dots+k_i\mypp,
\end{array}
\right.
\end{equation}
which corresponds to the measure
\begin{equation}\label{eq:tnu}
\nu(\diff x)=\sum_{j=1}^\infty \delta_{p_j}(\diff x)
=\sum_{i=0}^\infty k_i\, \delta_{k_{i+1}^{-1}}(\diff x)\myp.
\end{equation}
That is to say, the array of boxes is partitioned in blocks so that
$i$-th block contains $k_i$ boxes of frequencies $1/k_{i+1}$
\,($i=0,1,2,\dots$).

The heuristics underlying this example is as follows. A prototype
instance is a block of $k$ equal boxes each with frequency, say,
$q$. The mean number of singleton boxes within the block is a
single-wave function $k\myp tq\,\e^{-tq}$ which increases to its
maximum $k/\e$ at time $t=1/q$ and then goes down to~$0$. Now, the
idea is to combine a series of such blocks in order to guarantee a
suitable overlap of the waves produced by successive blocks. If the
sequence $(k_i)$ grows fast enough, then for each $i=0,1,2,\dots$
there exists a time instant (of order of $k_{i+1}$) when boxes
belonging to $i$-th block start to get occupied. After some time,
the mean number of singletons among these boxes is still relatively
large, say not less than $\log\log k_i$, but the expected number of
balls that fall in boxes of further blocks becomes large too, and
almost all these balls produce singleton boxes, since $k_{i+1}$ is
yet much larger (hence the frequencies are smaller). As time passes,
all boxes belonging to blocks $0,1,\dots,i$ are likely to contain
more than one ball each, while the balls hitting other blocks remain
sole representatives of their boxes.

To make this heuristic work, we choose
\begin{equation}\label{eq:k_i}
k_i:=2^{2^i},\qquad i=0,1,2,\dots,
\end{equation}
so that $k_{i+1}=k_i^2$ for all $i$. We wish to check that
$\Phi_1(t)$ goes to infinity but $\Phi_2(t)$ does not. Using
(\ref{eq:Phitr}) and (\ref{eq:tnu}) we have
\begin{align}\label{eq:exPhi1}
\Phi_1(t)=t\int_0^\infty x\,\e^{-tx}\,\nu(\diff x)=\sum_{i=0}^\infty
\frac{t\myp k_i}{k_{i+1}}\,
\e^{-t/k_{i+1}}&=:\sum_{i=0}^\infty A_i(t)\myp,\\
\label{eq:exPhi2} \Phi_2(t)=\frac{t^2}{2}\int_0^\infty x^2 \myp
\e^{-tx}\,\nu(\diff x)=\frac{1}{2}\sum_{i=0}^\infty
 \frac{t^2k_i}{k_{i+1}^2}\,\e^{-t/k_{i+1}}
&=:\frac{1}{2}\sum_{i=0}^\infty B_i(t)\myp.
\end{align}

As a function of $t$, each summand $A_i(t)$ in the sum
(\ref{eq:exPhi1}) increases up to the maximum value
$A_i(t_i^*)=k_i\myp \e^{-1}$ attained at $t_i^*=k_{i+1}$, and then
decreases to zero. Two consecutive summands, $A_i(t)$ and
$A_{i+1}(t)$, are equal at the point
$$
t_i':=\frac{k_{i+1}^2}{k_{i+1}-1}\log k_{i}\mypp,
$$
where their common value is
$$
A_i(t_i')=\frac{k_{i+1}}{k_{i+1}-1} \,k_i^{-1/(k_{i+1}-1)}\log
k_{i}\mypp.
$$
Using the elementary inequality $k^{-1/(k-1)}\ge \e^{-1}$
\,($k>1$)\myp, we note that
$$
A_i(t_i')\ge k_i^{-1/(k_{i}-1)}\log k_{i}\ge \e^{-1}\log k_{i}\mypp.
$$
Since $t_{i-1}'<t_i^*<t_i'$ \,$(i=1,2,\dots)$, it follows that for
all $t\in[t_{i-1}',t_i']$\myp,
$$
\Phi_1(t)\ge A_i(t)\ge \e^{-1} \log k_{i-1}\mypp,
$$
hence
$$
\liminf_{t\to\infty}\,\Phi_1(t)\ge \e^{-1}\liminf_{i\to\infty}
\,\log k_{i-1} =\infty\myp.
$$

Turning to $\Phi_2(t)$, note that the summand $B_i(t)$ in
(\ref{eq:exPhi2}) attains its maximum value at the point
$t=2t_i^{*}=2 k_{i+1}$ and $B_i(2t_i^*)=4\myp \e^{-2}\myp k_i$, so
$$
\Phi_2(2t_i^*)\ge B_i(2t_i^*)= 4\myp \e^{-2}\myp k_i\to
\infty\qquad(i\to\infty)\myp.
$$
On the other hand, on the sequence $t_j'':=3\myp k_{j+1}\log k_j$
one has
\begin{align*}
 B_i(t_j'') &= \frac{(t_j'')^2 k_i}{k_{i+1}^2}\,\e^{-t_j''/k_{i+1}}=
 \frac{9\myp k_{j+1}^2\log^2 k_j}{k_{i+1}^{3/2}}\,
 \exp\left(-\frac{3k_{j+1}\log k_j}{k_{i+1}}\right).
\end{align*}
Setting $x=k_{i+1}$ and $a=k_{j+1}\log k_j$, we note that the
function $x^{-3/2}\,\e^{-3a/x}$ increases for $0<x\le 2a$. Hence,
for all $i=0,1,\dots, j$\myp,
$$
 B_i(t_j'')\le B_j(t_j'')=\frac{9\log^2 k_j}{k_j^2}\mypp,
$$
and therefore
\begin{equation}\label{eq:j1}
 \sum_{i=0}^j B_i(t_j'')\le (j+1) B_j(t_j'')=
 \frac{9\myp(j+1)\log^2 k_j}{k_j^2}\mypp.
\end{equation}
For $i\ge j+1$, we have
$$
 B_i(t_j'')\le \frac{9\myp k_{j+1}^2\log^2 k_j}{k_ik_{i+1}}\le
 \frac{9\log^2 k_j}{k_i}\mypp,
$$
and since $k_i=2^{2^i}\ge 2^{4i}$ for $i\ge 4$, it follows
\begin{equation}\label{eq:j2}
\sum_{i=j+1}^\infty B_i(t_j'')\le 9\cdot 2^{2j}\sum_{i=j+1}^\infty
2^{-4i}=\frac{3}{5\cdot 2^{2j}}\mypp.
\end{equation}
Combining the estimates (\ref{eq:j1}) and (\ref{eq:j2}) yields
$\Phi_2(t_j'')\to 0$ as $j\to\infty$\myp.

Thus $\Phi_2(t)$ does not have a limit as $t\to\infty$, and moreover
$$
\liminf_{t\to\infty}\mypp\Phi_2(t)=0\myp,
\qquad\limsup_{t\to\infty}\mypp\Phi_2(t)=\infty\myp.
$$

Finally, it is easy to see directly that in this example the limit
in (\ref{eq:pjovertail}) does not exist. Indeed, along the
subsequence $j=k_0+k_1+\dots+k_i$, according to (\ref{eq:p_ex}) and
(\ref{eq:k_i}),
$$
\rho_j=k_{i+1}\left(\frac{k_{i+1}}{k_{i+2}}
 +\frac{k_{i+2}}{k_{i+3}}+\cdots\right)=1+O(k_{i+1}^{-1})\to
 1\qquad
(i\to\infty)\myp.
$$
On the other hand, for $j=k_0+k_1+\cdots+k_i+1$
we have
$$
\rho_j=k_{i+2}\left(\frac{k_{i+1}-1}{k_{i+2}}+
\frac{k_{i+2}}{k_{i+3}}+\cdots\right) \ge k_{i+1}-1\to \infty\qquad
(i\to\infty)\myp.
$$
\end{example}

Karlin \citep[page~384]{Karlin} gives an example of frequencies for
which $V(t)$ converges to $0$ along a sequence of values of $t$, and
converges to $\infty$ along another sequence; in that case
$\Phi_1(t)$ demonstrates the same type of behavior. Our Example
\ref{ex:Phi1notPhi2} exhibits a more exotic ``second order''
pathology: this time, $\Phi_1(t)\to \infty$ but $\Phi_2(t)$
oscillates between $0$ and $\infty$.

\subsection{Relationship between the various sufficient
conditions}\label{sec4.4}
First of all, note that condition (\ref{eq:(ii)}) in Theorem
\ref{suff1} does not imply condition (\ref{eq:(i)}). A
counterexample may be constructed by a slight modification of
Example~\ref{ex:3^k} as follows: define the frequencies $(p_j)$ by
setting $\nu(\diff x)=\sum_{i=1}^\infty i\delta_{\tilde p_i}$, where
$\tilde p_i:=i^{-1} 2^{-i}$, then $\liminf_{j\to\infty}
p_{j+k}/p_j=\liminf_{i\to\infty} \tilde p_{i+1}/\tilde p_i=1/2$ (so
that (\ref{eq:(ii)}) is satisfied), but for $(i+1)^{-1}\myp
2^{-i}\le x< i^{-1}\myp 2^{-i}$ we have $\Delta\nu(x)=0$\myp, hence
$\liminf_{x\downarrow 0}\Delta\nu(x)=0$ and (\ref{eq:(i)}) fails.

Further, it is easy to see that condition (\ref{eq:(ii)}) in Theorem
\ref{suff1} implies the set of equivalent conditions
(\ref{eq:nutoverx}), (\ref{eq:pjovertail}) in Theorem \ref{sctail},
but not the other way around. Indeed, if (\ref{eq:(ii)}) is
satisfied then for $\rho_j$ defined in (\ref{eq:rho}) we have
$$
\liminf_{j\to\infty} \rho_j \ge \liminf_{j\to\infty}\sum_{k=1}^M
\frac{p_{j+k}}{p_j}\ge M\cdot\frac{1}{2}\to\infty\qquad
(M\to\infty)\myp,
$$
and condition (\ref{eq:pjovertail}) follows. On the other hand, we
have seen that in Example \ref{ex:3^k} condition (\ref{eq:(ii)})
fails, while for $q^{j}\le x<q^{j-1}$ we have
$$
\frac{\tilde{\nu}[0,x]}{x}=\frac{1}{x}\sum_{i\ge j}^\infty i\myp
q^{i}\ge \frac{j}{q^{j-1}}\sum_{i\ge j}^\infty q^{i}=\frac{j\myp
q}{1-q}\to \infty\qquad(j\to\infty)\myp,
$$
and the condition (\ref{eq:nutoverx}) is valid.

As Example \ref{ex:Phi1notPhi2} shows, a converse to Theorem
\ref{sctail} is \emph{not valid}\/, unless under further assumptions
on the measure $\tilde{\nu}$ (cf.~\citep{Karlin,Dutko}). For
instance, if $\tilde{\nu}[0,x]$ varies regularly at zero, then
Karamata's Tauberian theorem (see
\citep[\S\myp1.7.2\myp]{Bingham-Goldie-Teugels} or \citep[\S\mypp
XIII.5\myp]{Feller2}) applied to (\ref{eq:Phitr-nutilde}) yields
$\tilde{\nu}[0,x]/x\sim c\,{\Phi_1(1/x)}$ as $x\downarrow 0$, so
that the convergence $\Phi_1(t)\to\infty$ as $t\to\infty$ does imply
the condition (\ref{eq:nutoverx}).

\begin{remark}
By Karamata's Tauberian theorem, the convergence
$$
\Phi_1(t)=t\int_0^\infty \e^{-tx}\,\tilde{\nu}(\diff x) \to
c\qquad(t\to\infty)
$$
is equivalent to $\tilde{\nu}[0,x]/x \to c$ as $x\downarrow 0$\mypp.
Interestingly, the implication may fail for $c=\infty$, as Example
\ref{ex:Phi1notPhi2} demonstrates.
\end{remark}

\section{Convergence to a finite limit}\label{sec5}
We will now investigate the situation where the variance $V(t)$ has
a finite limit as $t\to\infty$\myp, which is the central topic of
this work (see Theorem~\ref{thm:limit}). As already mentioned in
Section \ref{sec3.6}, the ``if'' part of Theorem \ref{thm:limit}
follows from Corollary \ref{cor:lagged}. So the main goal of this
section is to prove the ``only if'' part (i.e., the sufficiency of
the condition (\ref{pv})), but we will also give a streamlined proof
of the necessity.

\subsection{Criterion of convergence}
Recall that $D(\cdot)$ is a primitive function of
$\Delta\nu(\cdot)$, defined by (\ref{eq:D}).

\begin{lemma}\label{lm:Cesaro} In order that there exist a finite
limit
\begin{equation}\label{eq:var_limit}
\lim_{t\to\infty} V(t)=:v\myp,
\end{equation}
it is necessary and sufficient that
\begin{equation}\label{eq:Cesaro}
\lim_{x\downarrow 0} \frac{D(x)}{x}=v\myp.
\end{equation}
\end{lemma}
\proof Note that, according to (\ref{eq:vbar>0}), $v>0$. By the
representation (\ref{eq:V1}), we can rewrite (\ref{eq:var_limit}) as
\begin{equation}\label{eq:D_n}
\int_0^\infty \e^{-tx}\,\diff D(x)\sim \frac{v}{t}\qquad
(t\to\infty)\myp.
\end{equation}
By Karamata's Tauberian theorem (see
\citep[\S\myp1.7.2]{Bingham-Goldie-Teugels}, \citep[\S\mypp
XIII.5]{Feller2}), the relation (\ref{eq:D_n}) is equivalent to
$D(x)\sim vx$ as $x\downarrow 0$, which is the same as
(\ref{eq:Cesaro}).
\endproof

\subsection{Some implications of convergence}
\begin{lemma}\label{lm:alphabeta}
Suppose that the limit\/ \textup{(\ref{eq:Cesaro})}
exists\/\textup{,} and let\/ $\alpha,\,\beta>0$ be arbitrary
variables such that\/ $\alpha,\,\beta\downarrow0$ and\/
$(\alpha+\beta)/(\beta-\alpha)=O(1)$\myp. Then
%\begin{equation}\label{eq:D-D/beta-alpha}
$$
\lim_{\alpha,\myp\beta\downarrow0}
\frac{D(\beta)-D(\alpha)}{\beta-\alpha} = v\myp.
$$
%\end{equation}
\end{lemma}
\proof
Using (\ref{eq:Cesaro}), we have
$$
\frac{D(\beta)-D(\alpha)}{\beta-\alpha}=
\frac{v\beta\myp(1+o(1))-v\alpha\myp(1+o(1))}{\beta-\alpha}
=v+\frac{o(1)(\alpha+\beta)}{\beta-\alpha} \to v\myp,
$$
since the ratio $(\alpha+\beta)/(\beta-\alpha)$ is bounded.
\endproof

\begin{lemma}\label{lm:integer}
If the finite limit\/ {\rm (\ref{eq:var_limit})} exists then the
limiting value\/ $v$ must be a positive integer number\myp\textup{,}
$v=k\in\NN$\textup{,} and in this case
\begin{equation}\label{eq:lambda}
\lim_{x\downarrow 0}\frac{\lambda\{u\in{}]0,x]:\Delta\nu(u)\ne
k\}}{x}=0\myp,
\end{equation}
where\/ $\lambda\{\cdot\}$ denotes Lebesgue's measure on\/
$\RR_+$\mypp.
\end{lemma}
\proof By Lemma~\ref{lm:bound}, the function $\Delta\nu(u)$ is
uniformly bounded. By definition, it counts the number of
frequencies $p_j$ in the interval $]u/2,u]$, therefore
$\Delta\nu(u)$ is piecewise constant, with jumps at points $u=p_j$
and $u=2p_j$. Thus, for any given interval $]x/2,x]$ the total
number of such jumps is uniformly bounded by a constant, say
$M<\infty$.

Let $]\alpha,\beta[$ be the maximal open subinterval of $]x/2,x]$,
on which $\Delta\nu(\cdot)$ is constant. Clearly, its length
satisfies $\beta-\alpha\ge x/2\myp(M+1)$, thus
\begin{equation}\label{eq:M}
0\le\frac{\alpha+\beta}{\beta-\alpha}\le \frac{2x}{x/2\myp(M+1)}=
4\myp(M+1)\myp.
\end{equation}
Consider a closed interval $[\alpha_1,\beta_1]\subset{}
]\alpha,\beta[$ with $\alpha_1=(3\alpha+\beta)/4$,
$\beta_1=(3\beta+\alpha)/4$. Since $\alpha_1+\beta_1=\alpha+\beta$
and $\beta_1-\alpha_1=(\beta-\alpha)/2$, by the bound (\ref{eq:M})
Lemma~\ref{lm:alphabeta} applies to yield
\begin{equation}\label{eq:mean}
\frac{1}{\beta_1-\alpha_1}\int_{\alpha_1}^{\beta_1}
\Delta\nu(u)\,\diff
u=\frac{D(\beta_1)-D(\alpha_1)}{\beta_1-\alpha_1}\to
v\qquad(x\downarrow 0)\myp.
\end{equation}
But the function $\Delta\nu(\cdot)$ is constant on
$]\alpha,\beta[{}\supset[\alpha_1,\beta_1]$\myp, hence its sole
(integer) value must coincide with the asymptotic mean $v$ given by
(\ref{eq:mean}). In particular, $v$ must be integer, $v=k\in\NN$.

Along the same lines, one can show that for any $\varepsilon>0$ and
all small enough $x$, the function $\Delta\nu(\cdot)$ takes the
value $v=k$ on the interval $]x/2,x]$ everywhere except on a set of
Lebesgue's measure smaller than $\varepsilon x$. Thus, Lebesgue's
measure of the set $\{u\in{}]0,x]:\Delta\nu(u)\ne k\}$ is bounded by
$\varepsilon\sum_{i=1}^\infty 2^{-i+1}x=2\varepsilon x$, and since
$\varepsilon$ is arbitrary, (\ref{eq:lambda}) follows.
\endproof

\subsection{Lagged frequency ratio and
the proof of Theorem\/ \textup{\ref{thm:limit}}}

\begin{lemma}\label{th:pseudo}
If the limit\/ \textup{(\ref{eq:var_limit})} exists\/
\textup{(}hence $v=k\in\NN$ by Lemma\/
\textup{\ref{lm:integer}}\textup{),} then\/ \textup{(}cf.\
\textup{(\ref{pv})}\textup{)}
\begin{equation}\label{eq:pk/pkv}
\lim_{j\to\infty}\frac{p_{j+k}}{p_j}=\frac{1}{2}\mypp.
\end{equation}
\end{lemma}
\proof Without loss of generality, it suffices to consider $j\in\NN$
such that $2p_{j+k}\ne p_j$. Suppose first that $2p_{j+k}<p_j$. Then
for $x\in[2p_{j+k},p_j[$ we have $]x/2,x]\subset {}]p_{j+k},p_j[$
and hence $\Delta\nu(x)\le k-1$. Therefore,
\begin{equation}\label{eq:pk}
D(p_j) - D(2p_{j+k})= \int_{2p_{j+k}}^{p_j}\Delta\nu(u)\,\diff u \le
(k-1)(p_j-2p_{j+k})\myp.
\end{equation}
Using that $D(x)=k\myp x\left(1+o(1)\right)$ as $x\downarrow 0$ (see
Lemma \ref{lm:Cesaro}), from (\ref{eq:pk}) we deduce that
$\liminf_{j\to\infty}p_{j+k}/p_j\ge 1/2$, which, together with the
hypothesis $p_{j+k}/p_j<1/2$ (see above), implies (\ref{eq:pk/pkv}).

Likewise, if $p_j<2p_{j+k}$ then for $x\in[p_j,2p_{j+k}[$ we have
$]x/2,x]\supset [p_{j+k},p_j]$, hence $\Delta\nu(x)\ge k+1$ and
(cf.\ (\ref{eq:pk}))
\begin{align*}
D(2p_{j+k}) - D(p_j)= \int_{p_j}^{2p_{j+k}}\Delta\nu(u)\,\diff u \ge
(k+1)(2p_{j+k}-p_j)\myp.
\end{align*}
Similarly as before, this simplifies to
$\limsup_{j\to\infty}p_{j+k}/p_j\le 1/2$, and since we assumed that
$p_{j+k}/p_j<1/2$, (\ref{eq:pk/pkv}) follows. The proof is complete.
\endproof

Let us now show the converse of Lemma \ref{th:pseudo} (as mentioned
at the beginning of Section~\ref{sec5}, this also follows from
Corollary}~\ref{cor:lagged}).
\begin{lemma}\label{th:pseudo1}
Assume that the sequence $(p_j)$ satisfies the condition
\textup{(\ref{eq:pk/pkv})} for some $k\in\NN$. Then the limit
\textup{(\ref{eq:var_limit})} exists and $v=k$\myp.
\end{lemma}

\proof By additivity, it suffices to prove that for each subsequence
$p^{(i)}_j:=p_{i+k(j-1)}$ ($i=1,\dots, k$), its contribution to the
limit (\ref{eq:var_limit}) equals exactly~$1$. Thus the proof is
reduced to showing that if $(p_j)\in\RT_{1/2}$ then
\begin{equation}\label{eq:=1}
 V(t)=\sum_{j=1}^\infty \bigl(\e^{-tp_j}-\e^{-2tp_j}\bigr)\to1\qquad (t\to\infty)\myp.
\end{equation}

By the $\RT$-condition, $2p_{j+1}=p_j\myp(1+\gamma_j)$\myp, where
$\gamma_j\to0$ as $j\to\infty$. Hence, for any
$\varepsilon\in\myp]0,1/3]$ and all $j$ large enough we have
$|\gamma_j|\le\varepsilon$. In particular, $p_{j+2}/p_j\le
(1+\varepsilon)^{2}/4\le 4/9<1/2$\myp, which implies by Lemma
\ref{p-and-nu} that $\Delta\nu(x)\le2$ for small $x$. By
Lemma~\ref{lm:bound} and the estimate (\ref{PhiV}), it follows that
$\Phi_1(\cdot)$ is bounded. Returning to (\ref{eq:=1}), observe that
\begin{equation}\label{eq:gamma_k}
\sum_{j=j_0}^M \bigl(\e^{-tp_j}-\e^{-2tp_j}\bigr) =\sum_{j=j_0}^M
\e^{-tp_j}\bigl(1-\e^{-tp_j\gamma_j}\bigr)
-\e^{-2tp_{j_0}}+\e^{-2tp_{M+1}}.
\end{equation}
By the inequality $|1-\e^{-y}|\le |y|\mypp \e^{|y|}$\myp, the sum in
(\ref{eq:gamma_k}) is dominated by
\begin{align*}
\sum_{j=j_0}^M \e^{-tp_j(1-\varepsilon)}\myp tp_j\mypp \varepsilon
&\le\varepsilon\sum_{j=1}^\infty \e^{-tp_j(1-\varepsilon)}\myp tp_j
=\frac{\varepsilon}{1-\varepsilon}\,\Phi_1(t(1-\varepsilon)\bigr)
=O(\varepsilon).
\end{align*}
Passing to the limit in (\ref{eq:gamma_k}) as $M\to\infty$, we
obtain $V(t)=1+o(1)+O(\varepsilon)$ as $t\to\infty$, and since
$\varepsilon$ is arbitrarily small, we arrive at (\ref{eq:=1}).
\endproof

We are now able to complete the proof of our main Theorem
\ref{thm:limit} characterizing the case of converging variance.
Indeed, putting together Lemmas \ref{th:pseudo} and \ref{th:pseudo1}
yields the desired criterion for $V(t)\to v$. Appealing to Corollary
\ref{depoissmom} we conclude that the same condition applies to
$V_n\to v$.

\subsection{Link with Karlin's condition}\label{sec5Karlin}
In conclusion, let us recall that Karlin's sufficient condition for
$V(t)\to v$ \citep[Theorem 2]{Karlin} involves (i) the condition
$\limsup_{j\to\infty} p_{j+1}/p_j<1$ and (ii) an integral condition,
which in our notation reads
\begin{equation}\label{eq:Karlin0}
\lim_{x\to\infty} \frac{1}{x}\int_0^{x}\Delta\nu(1/y)\,\diff
y=v\myp,
\end{equation}
or, after an obvious change of variables,
\begin{equation}\label{eq:Karlin}
\lim_{x\downarrow 0}\,x\int_x^\infty\Delta\nu(u)\,u^{-2}\,\diff
u=v\myp.
\end{equation}

Throughout his paper, Karlin also postulates that the function
$\nu_{\rm c}(x)=\nu\mypp]x,\infty[$ is regularly varying at zero
(see \citep[pages 376--377]{Karlin}. As we shall see, this condition
is superfluous and may be omitted (in fact, Karlin's proof of his
Theorem~2 only requires the boundedness of $\Delta\nu(x)$, which
follows easily from condition (i)). Note that condition (i) itself
is not necessary for the convergence of $V(t)$: for instance, it
does not hold for a sequence $(p_j)$ obtained by merging several
geometric sequences with ratio $1/2$ into one.

Furthermore, application of condition (\ref{eq:Karlin0}) to the
geometric case (with ratio $q$) yields the following (cf.\
\citep[Example~6]{Karlin}
%page~385
containing an error). Let
$\log_{1/q}2=k+\delta$, where $k=[\log_{1/q}2]$ is the integer part
of $\log_{1/q}2$ and $\delta\in[0,1[$ is its fractional part. From
the definition of $\Delta\nu(\cdot)$ it follows that
\begin{align}
\notag\frac{1}{x}\int_0^x \Delta\nu(1/y)\,\diff y &=
\frac{1}{x}\int_0^x
\left([\log_{1/q}(2y)]-[\log_{1/q}y]\right)\diff y\\
\notag&= \frac{1}{x}\int_0^x \left([k+\delta +\log_{1/q}
y]-[\log_{1/q}y]\right)\diff y\\
\label{eq:k+}&= k+\frac{1}{x}\int_0^x
\left([\delta+\log_{1/q}y]-[\log_{1/q}y]\right)\diff y\myp.
\end{align}
If $\delta=0$, the integral in (\ref{eq:k+}) vanishes and condition
(\ref{eq:Karlin0}) yields $v=k$. However, if $0<\delta<1$ then
(\ref{eq:k+}) does not have a limit as $x\to\infty$, since for
$x=q^{-j}$ the integral term amounts to
$$
q^j\sum_{i=1}^j q^{-i}(1-q^\delta)\to
\frac{1-q^\delta}{1-q}\qquad(j\to\infty)\myp,
$$
whereas for $x=q^{-j-1+\delta}$ it reads
$$
q^{j+1-\delta}\sum_{i=1}^j q^{-i}(1-q^\delta)\to
q^{1-\delta}\,\frac{1-q^\delta}{1-q}\qquad(j\to\infty)\myp.
$$
As a result, condition (\ref{eq:Karlin0}) is satisfied if and only
if $\log_{1/q}2=k\in\NN$\myp, or equivalently $q=2^{-1/k}$. Our
Theorem~\ref{thm:limit} gives the same result, so (\ref{eq:Karlin0})
proves to yield a correct answer in the whole range of the geometric
case.

This observation brings up the question about the exact relationship
between Karlin's condition (\ref{eq:Karlin0}) (or (\ref{eq:Karlin}))
and our criterion (\ref{eq:Cesaro}). Surprisingly enough, we can
demonstrate the following.
\begin{theorem}\label{th:Karlin}
Condition\/ \textup{(\ref{eq:Karlin})} is equivalent to\/
\textup{(\ref{eq:Cesaro})}\textup{,} and hence the former is
necessary and sufficient in order that\/ $V(t)\to v$ as\/
$t\to\infty$.
\end{theorem}
\proof Suppose condition (\ref{eq:Cesaro}) holds. Using the notation
$D(x)$ (see (\ref{eq:D})) and integrating by parts, we get
\begin{multline*}
 x\int_x^\infty \Delta\nu(u)\,\frac{\diff u}{u^2}=x\int_x^\infty
 u^{-2}\,\diff D(u)=-\frac{D(x)}{x}+2x\int_x^\infty D(u)\,u^{-3}\,\diff u\\
 =
-\frac{D(x)}{x}+2\int_1^\infty \frac{D(xs)}{xs}\,s^{-2}\,\diff s \to
-v+2v\int_1^\infty s^{-2}\,\diff s=v\quad (x\downarrow 0)\myp,
\end{multline*}
where we used that the function $D(u)/u$ is bounded on $]0,\infty[$
(in particular, the dominated convergence theorem can be applied).
Hence, (\ref{eq:Karlin}) follows.

On the other hand, condition (\ref{eq:Karlin}) amounts to
\begin{equation}\label{eq:xG}
\lim_{x\downarrow 0} x\myp G(x)=v\myp,\qquad G(x):=\int_x^\infty
\Delta\nu(u)\,u^{-2}\,\diff u\myp.
\end{equation}
Again integrating by parts, we obtain
\begin{align*}
\frac{1}{x}\int_0^x \Delta\nu(u)\,\diff u &=-\frac{1}{x}\int_0^x
u^2\,\diff G(u)=-x\myp G(x)+\frac{2}{x}\int_0^x u\mypp G(u)\,\diff u\\
&=-x\myp G(x)+2\int_0^1 xs\,G(xs)\,\diff s\to -v+2v=v\quad
(x\downarrow 0)\myp,
\end{align*}
where we may use dominated convergence because the function $u\mypp
G(u)$ is bounded on $]0,1]$ due to (\ref{eq:xG}). Thus, condition
(\ref{eq:Karlin}) implies (\ref{eq:Cesaro}), and the proof is
complete.
\endproof

\begin{remark}
The statement of Theorem \ref{th:Karlin} is a particular case of a
general Karamata theorem (see
\citep[\S\mypp1.6.3]{Bingham-Goldie-Teugels},
\citep[\S\mypp{}VIII.9]{Feller2}), according to which the limiting
relation (\ref{eq:Cesaro}) is equivalent to either of the limits
\begin{align*}
\lim_{x\downarrow0} \, x^{\sigma-1}\!\int_x^\infty
\Delta\nu(u)\,u^{-\sigma}\,\diff
u&=\frac{v}{\sigma-1}\qquad(\sigma>1)\myp,\\
\lim_{x\downarrow0} \, x^{\sigma-1}\!\int_0^x
\Delta\nu(u)\,u^{-\sigma}\,\diff
u&=\frac{v}{1-\sigma}\qquad(\sigma<1)\myp.
\end{align*}
(Note that (\ref{eq:Cesaro}) itself is contained in the second
formula with $\sigma =0$.) That is to say, our condition
(\ref{eq:Cesaro}) may be included in a parametric family of mutually
equivalent criteria, set in terms of rescaled integrals of the
function $\Delta\nu(\cdot)$ against polynomial weights (the
canonical criterion (\ref{eq:Cesaro}) being apparently the
simplest). We have given a direct proof of Theorem \ref{th:Karlin}
because of the historic interest of Karlin's
condition~(\ref{eq:Karlin0}).
\end{remark}

\section*{Acknowledgment}
Main part of this research was carried out 
during L.~V.~Bogachev's visit
to the University of Utrecht in June 2006, made possible due to a
grant under the Global Exchange Programme (GEP) of the World
Universities Network (WUN), awarded at the University of Leeds.
Hospitality of the hosts at Utrecht is much appreciated.

\end{document}